\documentclass[10pt]{article}
\usepackage{amssymb,amsmath,amsthm}
\usepackage[numbers,sort&compress]{natbib}
\usepackage{indentfirst}
\usepackage{exscale}
\usepackage{relsize}
\usepackage{subfigure}
\usepackage{graphicx}
\usepackage{epstopdf}
\usepackage{xcolor}

\newcommand{\R}{\mathbb{R}}
\theoremstyle{definition}

\theoremstyle{remark}

\numberwithin{equation}{section}
\allowdisplaybreaks

 \UseRawInputEncoding

\begin{document}
\title{\Large\bf{ Ground state sign-changing solutions for Kirchhoff-type equations with logarithmic nonlinearity on locally finite graphs}}
\date{}
\author {Xin Ou$^1$, \ Xingyong Zhang$^{1,2}$\footnote{Corresponding author, E-mail address: zhangxingyong1@163.com}\\
{\footnotesize $^1$Faculty of Science, Kunming University of Science and Technology,}\\
 {\footnotesize Kunming, Yunnan, 650500, P.R. China.}\\
{\footnotesize $^{2}$Research Center for Mathematics and Interdisciplinary Sciences, Kunming University of Science and Technology,}\\
 {\footnotesize Kunming, Yunnan, 650500, P.R. China.}\\
}

 \date{}
 \maketitle

 \begin{center}
 \begin{minipage}{15cm}
 \small  {\bf Abstract:} We obtain the existence results of ground state sign-changing solutions and ground state solutions for a class of Kirchhoff-type equations with logarithmic nonlinearity on a locally finite graph $G=(V,E)$, and obtain the sign-changing ground state energy is  larger than twice of the ground state energy. The method we used is a direct non-Nehari manifold method in  [X.H. Tang, B.T. Cheng. J. Differ. Equations. 261(2016), 2384-2402.]

 \end{minipage}
 \end{center}
  \allowdisplaybreaks
 \vskip2mm
 {\section{Introduction }}
\setcounter{equation}{0}

\par
Let $G=(V, E)$ be a graph, where $V,E$ denote the vertex set and the edge set, respectively. A graph is said to be {\it locally finite} if for any $x\in V$, there exist only finite $y\in V$ such that $xy\in E$. A graph is {\it connected} if any two vertices $x$ and $y$ can be connected via finite edges. For any $x,y\in V$ with $xy\in E$, assume that its weight satisfies $\omega_{xy}=\omega_{yx}>0$. For any $x\in V$, define its degree as $deg(x)=\sum_{y\thicksim x}\omega_{xy}$, where we write $y\thicksim x$ if $xy\in E$. The measure $\mu:V\rightarrow \R^+$ is a finite positive function.
The {\it distance $d(x,y)$} of two vertices $x,y\in V$ is defined by the minimal number of edges which connect $x$ and $y$. We call that $\Omega\subset V$ is a {\it bounded domain} in $V$, if the distance $d(x,y)$ is uniformly bounded from above for any $x,y\in\Omega$. We denote the boundary of $\Omega$ by $\partial\Omega$ which is defined as
$$
\partial\Omega:=\{y\in V,\;y\notin\Omega:\exists\;x\in\Omega  \;\;\text{such that}\; xy \in E\}
$$
and the interior of $\Omega$ by $\Omega^\circ$. Then it is easy to see that $\Omega=\Omega^\circ$ which is different from that in the Euclidean  setting.
\par
For any function $u: \Omega\rightarrow \mathbb{R}$, define
\begin{eqnarray}
\label{eq1}
\Delta u(x)=\frac{1}{\mu(x)}\sum\limits_{y\thicksim x}w_{xy}(u(y)-u(x)).
\end{eqnarray}
The corresponding gradient is defined as
\begin{eqnarray}
\label{eq2}
\Gamma(u_1,u_2)(x)=\frac{1}{2\mu(x)}\sum\limits_{y\thicksim x}w_{xy}(u_1(y)-u_1(x))(u_2(y)-u_2(x)):=\nabla u_1 \cdot\nabla u_2.
\end{eqnarray}
Denote $\Gamma(u)=\Gamma(u,u)$ and  the length of the gradient is defined by
\begin{eqnarray}
\label{eq3}
|\nabla u|(x)=\sqrt{\Gamma(u)(x)}=\left(\frac{1}{2\mu(x)}\sum\limits_{y\thicksim x}w_{xy}(u(y)-u(x))^2\right)^{\frac{1}{2}}.
\end{eqnarray}
For any function $u:\Omega\rightarrow\mathbb{R}$, define
\begin{eqnarray}
\label{eq4}
\int_\Omega u(x) d\mu=\sum\limits_{x\in \Omega}\mu(x)u(x).
\end{eqnarray}
In the distributional sense, $\Delta u$ can be expressed  as follows: for any $u\in C_c(\Omega)$,
\begin{eqnarray*}
\int_{\Omega\cup\partial\Omega}(\Delta u)\phi d\mu=-\int_{\Omega\cup\partial\Omega}\Gamma(u,\phi)d\mu,
\end{eqnarray*}
where $C_c(\Omega)$ is the set of all real valued continuous functions with compact support. The above knowledge and statements can be seen in
 \cite{Grigor 2017} and \cite{Yamabe 2016}.

\par
The energy relationship of solutions to Kirchhoff-type problems has always been of great interests (see \cite{Cheng 2020, Feng 2021, Ji 2023, Shuai 2015, Tang 2016}). Next, we only recall two works which are directly related to our results in this paper. In \cite{Tang 2016}, Tang-Cheng investigated the equation (\ref{eq5}) with $\alpha=1$, that is, the following Kirchhoff-type equation:
\begin{eqnarray}
\label{eq6}
  \begin{cases}
  -(a+b \int_{\Omega} |\nabla u|^2dx)\Delta u=f(x,u), & \text {in} \; \Omega,\\
  u=0, & \text {on} \; \partial\Omega,
  \end{cases}
\end{eqnarray}
They assumed that $f$ has a super-cubic growth  and a weaker monotonicity condition than these in \cite{Shuai 2015}. They obtained  that equation (\ref{eq6}) has a ground state solution and a ground state sign-changing solution with two precisely nodal domains  via constraint variational method, topological degree theory and some energy estimate inequalities. They also found that the energy of  ground state sign-changing solution is strictly larger than twice that of the ground state solutions of Nehari-type.  In particular, they emphasise that $f\in C(\R,\R)$ rather than $f\in C^1(\R,\R)$ in \cite{Shuai 2015}. But in fact the condition $f\in C^1$ is still used in their proof. However, their approach called as direct non-Nehari manifold method, is more direct than the one in \cite{Shuai 2015}. The real sense of the weakening the condition $f\in C^1$ to $f\in C$ is given in \cite{Cheng 2020}, where Cheng-Chen-Zhang researched the following Kirchhoff-type Laplacian equation  with Dirichlet boundary value condition:
\begin{eqnarray}
\label{eq5}
  \begin{cases}
  -M(\int_{\Omega} |\nabla u|^2dx)\Delta u=f(x,u), & \text {in} \; \Omega,\\
  u=0, & \text {on} \; \partial\Omega,
  \end{cases}
\end{eqnarray}
where $\Omega\subset \R^N$ $(N=1,2,3)$ is a bounded domain and $M$ is a Kirchhoff-type function in the form of $M(t)=a+bt^\alpha$ for $t\in\mathbb{R}^+,a,\alpha>0,b\geq0$. When $f$satisfies the super $2\alpha+1$ times growth and a weaker assumption than the usual Nehari-type monotonicity condition, they combined the method in \cite{Tang 2016} with the variant of Miranda theorem to obtain the equation (\ref{eq5}) has a least energy nodal solution. They also obtained the result of the twofold energy relationship.
\par
The existence of solutions for the Kirchhoff-type equations with the general form $f$ becoming logarithmic nonlinearity in the Euclidean setting is a hot topic in recent years and we refer readers to \cite{Bu 2023, Ding 2021, Gao 2023, Liang 2020, Shao 2021, Wen 2019, Li 2020}. Next, we only recall three works which are directly related to our results in this paper. In \cite{Wen 2019}, Wen-Tang-Chen used the method in \cite{Tang 2016} and generalized Tang-Cheng's results to  the following Kirchhoff equation with  logarithmic nonlinearity:
\begin{eqnarray}
\label{r3}
  \begin{cases}
  -(a+b \int_{\Omega} |\nabla u|^2dx)\Delta u+W(x)u=|u|^{p-2}u\ln u^2, & \text {in} \; \Omega,\\
  u=0, & \text {on} \; \partial\Omega
  \end{cases}
\end{eqnarray}
where $a,b>0$, $4<p<2^*$ and $\Omega\subset \R^3$ is a bounded domain. Furthermore, in \cite{Li 2020}, Li-Wang-Zhang also used the method in \cite{Tang 2016} and generalized Wen-Tang-Chen's results in \cite{Wen 2019} to  the following $p$-Laplacian Kirchhoff-type equation with logarithmic nonlinearity:
\begin{eqnarray*}
  \begin{cases}
  -(a+b\int_{\Omega} |\nabla u|^pdx)\Delta_p u=|u|^{q-2}u\ln u^2, & \text {in} \; \Omega,\\
  u=0, & \text {on} \; \partial\Omega.
  \end{cases}
\end{eqnarray*}
where $a,b>0$ are constant, $4\leq 2p<q<p^*$, $\Omega\subset \mathbb{R}^N$ is a bounded domain and $N>p$. In \cite{Gao 2023}, Gao-Jiang-Liu-Wei studied the following logarithmic Kirchhoff-type problem:
\begin{eqnarray}
\label{r1}
  \left(a+b \int_{\mathbb{R}^3} \big[|\nabla u|^2+V(x)u^2\big]dx\right)\left[-\Delta u+V(x)u\right]=|u|^{p-1}u\log |u|,\;\;x\in\mathbb{R}^3,
\end{eqnarray}
where $a,b>0$ and $1<p<3$. When the potential function $V(x)$ is coercive and bounded from below, using the variational methods and the maximum principle, they obtained the equation (\ref{r1}) has only trivial solution for all large enough $b$ and two positive solutions for all sufficiently small $b$.
\par
When $a=1$ and $b=0$ in equation (\ref{eq5}), the Kirchhoff-type equation becomes the second order elliptic partial differential equation without nonlocal term. The elliptic partial differential equations on graphs have many applications in image processing, data processing, machine learning, electric network, etc. (see \cite{Chung2005, Cheung 2018, Ta 2011, Elmoataz 2012}).  In recent years, the existence of solutions for the equations on locally finite graphs have attracted some attention. We refer readers to \cite{Grigor 2017, Zhang 2018, Han 2021, Chang 2023, Changxiaojun, Yamabe 2016}. Next, we only recall three works which are directly related to our results in this paper. In \cite{Grigor 2017}, Grigor'yan-Lin-Yang investigated the following Yamabe equation on a locally finite graph $G = (V, E)$:
\begin{eqnarray}
\label{r2}
 -\Delta u+h(x)u=f(x,u), \;\;\;\;   x \in  V,
\end{eqnarray}
where $h:V\rightarrow \mathbb{R}$ and $f:V\times \mathbb{R}\rightarrow \mathbb{R}$. They assumed that the measure $\mu(x)\geq\mu_{\min}=\min_{x\in V}\mu(x)>0$ for all $x\in V$. Let $h$ and $f$ satisfy the following conditions:
\par
{\it
$(H_1)$ there exists a constant $h_0>0$ such that $h(x)\geq h_0$ for all $x\in V$;
\par
$(H_2)$  $\frac{1}{h} \in L^{1}(V)$;
\par
$(F_1)$  $f(x,s)$ is continuous in $s$, $f(x,0)=0$, and for any fixed $M>0$, there exists a constant $A_M$ such that $\max_{s\in[0, M]}f(x,s)\leq A_M$ for all $x\in V$;
\par
$(F_2)$  $\limsup_{s\rightarrow 0^+}\frac{2F(x,s)}{s^2}<\lambda_1=\inf_{\int_V u^2d\mu=1}\int_V (|\nabla u|^2+hu^2)d\mu$;
\par
$(F_3)$  there exists a constant $\theta >2$ and $s>0$ such that
$$
  0<\theta F(x,s)=\theta\int_0^s f(x,t)dt\leq sf(x,s),\;\;\forall x\in V.
$$}
Then they obtained that equation (\ref{r2}) has a strictly positive solution by mountain pass theorem.
 \par
 In  \cite{Chang 2023}, Chang-Wang-Yan researched (\ref{r2}) with $f(x,u)=u\log {u^2}$, when $h$ satisfies the following conditions:
\par
{\it
$(H'_1)$ $h:V\rightarrow \mathbb{R}$ satisfies $\inf_{x\in V}h(x)\geq h_0$ for some constant $h_0\in(-1,0)$;
\par
$(H_3)$ for every $M>0$, the volume of the set $D_M$ is finite, namely,
$$
  Vol(D_M)=\sum_{x\in D_M}\mu (x)<\infty,
$$
where $D_M=\{x\in V:h(x)\leq M\}$.}\\
They obtained a ground state solution of (\ref{r2}) by Nehari manifold method. When $h$ satisfies the $(H'_1)$ and
\par
{\it
$(H'_3)$ there exists $M_0>0$ such that $1/h(x)\in L^1(V\backslash D_{M_0})$, where $D_{M_0}=\{x\in V:h(x)\leq M_0\}$ and the volume of $D_{M_0}$ is finite. }\\
They also obtained a ground state solution of (\ref{r2}) by the mountain pass theorem.
 \par
In \cite{Changxiaojun}, Chang-R\u{a}dulescu-Wang-Yan studied the (\ref{r2}) with $h(x)=\lambda a(x)$, where $\lambda$ is a positive constant. When $a$ satisfies the following condition:
\par
{\it
$(H_4)$ $a(x)\geq0$ and the potential well $\Omega=\{x\in V:a(x)=0\}$ is a non-empty, connected and domain in $V$.
\par
$(H_5)$ there exist $M>0$ such that the volume of the set $D_M$ is finite, namely,
$$
  Vol(D_M)=\sum_{x\in D_M}\mu (x)<\infty,
$$
where $D_M=\{x\in V:a(x)< M\}$.}\\
They proved (\ref{r2}) has a ground state sign-changing solution by variational techniques and the Nehari manifold method.
\par
In \cite{Ji 2023}, Pan-Ji investigated the following Kirchhoff equations on locally finite graphs:
\begin{eqnarray}
\label{r4}
-(a+b\int_V |\nabla u|^2d\mu)\Delta u+c(x)u=f(u),\;\; x\in V,
\end{eqnarray}
where $V$ is a locally finite graph and $a,b$ are positive constants. They assumed that $c$ satisfies $(H_1)$ and
\par
{\it
$(H_5)$ $c(x)\to +\infty$, as $dist(x,x_0)\to +\infty$ for some fixed $x_0\in V$.}\\
Moreover, they let $f$ satisfy the following conditions:
\par
{\it
$(F_1)$ $f(s)=o(|s|)$ as $s\to 0$;
\par
$(F_2)$ For $4<p<+\infty$, $\lim_{|s|\to+\infty}\frac{f(s)}{|s|^{p-1}}=0$;
\par
$(F_3)$ $\lim_{|s|\to+\infty}\frac{F(s)}{|s|^{4}}=+\infty$, where $F(s)=\int_0^s f(t)dt$;
\par
$(F_4)$ $\frac{f(s)}{|s|^{3}}$ is an increasing function of $s\in\mathbb{R}\setminus\{0\}$.}\\
Then they used the constrained variational method to prove the equation (\ref{r4}) has a ground state sign-changing solution and the energy of ground state sign-changing solution is strictly larger than twice that of the ground state solution.
\par
In this paper, inspired by \cite{Cheng 2020, Chang 2023, Yamabe 2016, Tang 2016, Ji 2023, Wen 2019, Li 2020},  we are concerned with the following Kirchhoff-type equation with logarithmic nonlinearity and Dirichlet boundary value condition on a locally finite graph $G=(V,E)$:
\begin{eqnarray}
\label{eq1}
 \begin{cases}
  -\left(a+b\int_{\overline{\Omega}}|\nabla u|^2d\mu\right)\Delta u+\lambda g(x)u^{\frac{2k}{m}-1}=Q(x)|u|^{p-2}u\ln |u|^r,& \text {in} \; \Omega^\circ,\\
  u=0,&\text {on} \; \partial\Omega,\\
   \end{cases}
\end{eqnarray}
where $\overline{\Omega} \subset V $ is a bounded domain, $\overline{\Omega}:=\Omega^\circ\cup\partial \Omega=\Omega\cup\partial \Omega$, $a>0$, $b\geq0$, $p>4$, $\lambda\ge0$, $r\geq 1$, $m,k\in \mathbb{N}$ satisfying $1<\frac{2k}{m}\leq p$ and $Q,g\in\mathcal{C} (\Omega^\circ, (0,\infty))$. We establish the existence results of ground state solutions and ground state sign-changing solutions for problem (\ref{eq1}) on the locally finite graph by using the method in \cite{Tang 2016} and \cite{Cheng 2020}, and obtain the sign-changing ground state energy is  larger than twice of the ground state energy. The main difficulties for this problem are that
$$
 I(u) \not =   I(u^+)+I(u^-) \;\;\text{and}\;\; \langle I'(u),u^\pm\rangle
\not=  \langle I'(u^\pm),u^\pm\rangle
$$
because of the special definition of  $|\nabla u|$ (see (\ref{eq3})), where
$$
u^+(x):=\max\{u(x),0\}\;\;\text{and}\;\; u^-(x):=\min\{u(x),0\}
$$
and $I$ is defined by (\ref{eq9}) below.
The phenomenon is different from the usual second order elliptic PDE like (\ref{eq6}) with $b=0$ in the Euclidean  setting.
However, by some careful computation,  we find two interesting relationships which are similar to the Kirchhoff-type PDE like (\ref{eq6}) and (\ref{r3}) in the Euclidean  setting (in some sense) (see Proposition 2.2 below) so that we can use the method in \cite{Tang 2016}  and \cite{Cheng 2020} to prove our results.
To be precise, we obtain the following result.
\par
\vskip2mm
\noindent
{\bf Definition 1.1.}  {\it If $u\in \mathcal{H}^{1,2}_0(\Omega)$ is a solution of equation (\ref{eq1}) and $u^\pm\not\equiv 0$, then $u$ is a sign-changing solution of (\ref{eq1}).}

\vskip2mm
\noindent
{\bf Theorem 1.1.} {\it Assume that $G=(V,E)$ is a locally finite graph,  $\Omega^\circ\not=\emptyset$
 and $\partial \Omega\not=\emptyset$. Then for each $\lambda\in[0,+\infty)$, problem} (\ref{eq1}) {\it has a sign-changing solution $\tilde{u}\in \mathcal{M}$ such that $I(\tilde{u})=\inf_{\mathcal{M}}I:=m$ and a constant-sign solution $\bar{u}\in \mathcal{N}$ such that $I(\bar{u})=\inf_{\mathcal{N}}I:=c$, where
$$
\mathcal{M}=\left\{
u\in \mathcal{H}^{1,2}_0(\Omega),\;u^{\pm}\neq 0,\;\;\langle I'(u),u^{+}\rangle=0\;\;and\;\;\langle I'(u),u^{-}\rangle=0
\right\}
$$
and
$$
\mathcal{N}=\{u\in \mathcal{H}^{1,2}_0(\Omega),\;u\neq 0\;and\; \langle I'(u),u\rangle=0\}.
$$
Moreover, $m\geq 2c$.}

\vskip2mm
\noindent
{\bf Remark 1.1.}   We generalize the results with $r=2$ in \cite{Wen 2019} to the locally finite graphs setting with more general nonlinearity. However, we do not obtain the result that the ground state sign-changing solution has two precisely nodal domains. The main reason is that the special definition of $G=(V,E)$ cause that a point $x_1\in V$  satisfying $\tilde{u}(x_1)>0$ is possible to connect directly to another point  $x_2\in V$  satisfying $\tilde{u}(x_2)<0$. This is different from the case in the Euclidean setting where  $\tilde{u}$ continuously varies from the part with  $\tilde{u}>0$ to the part with  $\tilde{u}<0$ and it has to pass through the part with  $\tilde{u}=0$. Maybe, it is possible to solve this problem for some graphs with some special structures, but currently, we have no idea for this.
\vskip2mm
\noindent
{\bf Remark 1.2.} The logarithmic nonlinearity in (\ref{eq1}) is more general than that in (\ref{r3}), that is, when $Q(x)\equiv 1$ for all $x\in V$ and $r=2$,  the logarithmic nonlinearity in (\ref{eq1}) reduces to that in (\ref{r3}). Moreover, (\ref{eq1}) is allowed to have more general power law, that is, $u^{\frac{2k}{m}-1}$ reduces to $u$ in (\ref{r2}) if $k=m=1$. Hence, in the setting of locally finite graphs, we allow more types of equations to have ground state sign-changing solution. For example, let $a=1$, $b=0$, $k=2$, $m=1$, $p=5$ and $r=1$. Then the following equation satisfies our Theorem 1.1:
\begin{eqnarray*}
\label{eqq1}
 \begin{cases}
  -\Delta u+\lambda g(x)u^3=Q(x)|u|^{3}u\ln |u|,& \text {in} \; \Omega^\circ,\\
  u=0,&\text {on} \; \partial\Omega,\\
   \end{cases}
\end{eqnarray*}
where $\lambda\ge 0$ and $g,Q:\Omega^\circ\to (0,+\infty)$.

\vskip2mm
\noindent
{\bf Remark 1.3.} Theorem 1.1 is different from that in \cite{Changxiaojun} where the sign-changing ground state solutions of  (\ref{r2}) with $f(x,u)=u\log {u^2}$ and $h(x)=\lambda a(x)$, where $\lambda$ is a parameter, are also investigated by using the similar approach, since we allows a different logarithmic nonlinearity and power law in  $({\ref{eq1}})$ as stated in Remark 1.2 and we consider the more general Kirchhoff type which contains a nonlocal term. Moreover, $({\ref{eq1}})$ is  defined on a bounded domain $\Omega\subset V$ rather than on the whole locally finite graph $G=(V,E)$ like (\ref{r2}), which, together with the Kirchhoff term, causes that the establishment of the relationship among $I(u)$, $I(u^+)$ and $I(u^-)$ are more complex and difficult than (\ref{r2}) in \cite{Changxiaojun}. Another difference is that the parameter $\lambda$ in Theorem 1.1 is allowed to belong to $[0,+\infty)$ rather than $[\lambda_0,+\infty)$ for some $\lambda_0>0$ which are needed in \cite{Changxiaojun}.

\vskip2mm
\noindent
{\bf Remark 1.4.}  We would like to address that Theorem 1.1 is also different from that in \cite{Ji 2023} where the sign-changing ground state solutions of  (\ref{r4}) is investigated. On one hand, our problem is  defined on a bounded domain $\Omega\subset V$ as described in Remark 1.3 and  we allow to have more general power law  $u^{\frac{2k}{m}-1}$ as described in Remark 1.2. On the other hand, our nonlinearity
$Q(x)|u|^{p-2}u\ln |u|^r$ does not necessarily satisfy $(f_4)$ in \cite{Ji 2023} even if $Q(x)\equiv 1$ for all $x\in \Omega$. For example, let $p,r=(7,5)$, $(6,6)$, $(9,2)$ and $(8,3)$ respectively. Then $\frac{|s|^{p-2}s\ln |s|^r}{|s|^3}$ is not increasing for $s\in \R\backslash\{0\}$ (see Figure 1).
\begin{figure*}[t]
	\centering
	\includegraphics[width=0.8\textwidth]{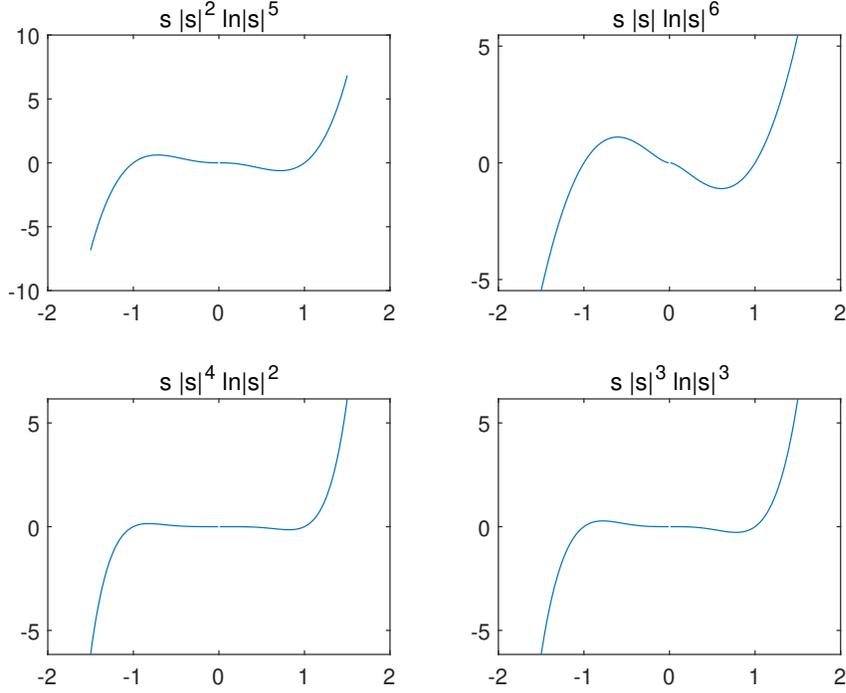}
	\vskip 0.2cm
	\caption{The diagram of function $\frac{|s|^{p-2}s\ln |s|^r}{|s|^3}$ with different values}
	\label{Figure1}
\end{figure*}

\vskip2mm
{\section{Variational setting}}
\setcounter{equation}{0}

Define
\begin{eqnarray}
\label{eq8}
\mathcal{H}_{0}^{1,2}(\Omega)=\left\{u:\Omega \rightarrow \mathbb{R}:\int_{\Omega\cup \partial \Omega}|\nabla u(x)|^2d\mu<\infty,\;u|_{\partial \Omega}=0\right\},
\end{eqnarray}
which is endowed with the norm
\begin{eqnarray*}
\|u\|=\left(\int_{\Omega\cup \partial \Omega}|\nabla u(x)|^2d\mu\right)^\frac{1}{2}.
\end{eqnarray*}
Then $\mathcal{H}_{0}^{1,2}(\Omega)$ is a Banach space and of finite dimension. Let $2\leq r<+\infty$. Define
$$
L^r(\Omega)=\left\{u:\Omega\to\mathbb{R}:\int_\Omega|u(x)|^rd\mu<\infty\right\}
$$
with the norm
\begin{eqnarray*}
\|u\|_{L^r(\Omega)}=\left(\int_\Omega|u(x)|^rd\mu\right)^\frac{1}{r}.
\end{eqnarray*}
When $r=+\infty$, we define
$$
L^\infty(\Omega)=\left\{u:\Omega\to\mathbb{R}:\sup_{x\in\Omega}|u(x)|<\infty\right\}
$$
with
$$
\|u\|_{L^\infty(\Omega)}=\max_{x\in\Omega}|u(x)|.
$$
{\bf Lemma 2.1.} \cite{Yamabe 2016}   {\it Let $G=(V,E)$ be a finite graph and $\Omega$ be a bounded domain satisfying $\Omega^\circ\not=\emptyset$. Then $\mathcal{H}^{1,2}_0(\Omega)\hookrightarrow L^q(\Omega)$ for all $1\leq q\leq+\infty$. Especially, if $1\leq q\leq +\infty$, then for all $\psi\in \mathcal{H}^{1,2}_0(\Omega)$,
\begin{eqnarray}
\label{Eq11}
\|\psi\|_{L^q(\Omega)}\leq K_q\|\psi\|,
\end{eqnarray}
where
$$
K_q=\frac{\left(\sum_{x\in V}\mu(x)\right)^{\frac{1}{q}}}{\mu_{\min}^{\frac{1}{2}}}, \ \mu_{\min}=\min_{x\in\Omega}\mu(x).
$$
In addition, $\mathcal{H}^{1,2}_0(\Omega)$ is pre-compact, namely, if $\{\psi_k\}$ is bounded in $\mathcal{H}^{1,2}_0(\Omega)$, then up to a subsequence, there exists some $\psi\in \mathcal{H}^{1,2}_0(\Omega)$ such that $\psi_k\rightarrow \psi$ in $\mathcal{H}^{1,2}_0(\Omega)$.}
\vskip2mm
\par
Define the functional $I:\mathcal{H}_{0}^{1,2}(\Omega)\to\R$ as
\begin{eqnarray}
\label{eq9}
    I(u)
&=& \frac{a}{2}\int_{\Omega\cup \partial \Omega} |\nabla u|^2d\mu
       +\frac{b}{4}\left(\int_{\Omega\cup \partial \Omega} |\nabla u|^2d\mu\right)^2
       +\frac{m}{2k}\int_\Omega\lambda g(x)u^{\frac{2k}{m}}d\mu\nonumber\\
& & +\frac{r}{p^2}\int_\Omega Q(x)|u|^pd\mu-\frac{1}{p}\int_\Omega Q(x)|u|^{p}|\ln |u|^r d\mu.
\end{eqnarray}
Note that
\begin{eqnarray*}
\label{eq10}
\lim_{t\rightarrow 0}\frac{Q(x)t^{p-1}\ln |t|^r}{t}=0&\text {and}&\lim_{t\rightarrow \infty}\frac{Q(x)t^{p-1}\ln |t|^r}{t^{l-1}}=0
\end{eqnarray*}
for all $x\in \Omega$, where $l\in(p,\infty)$. Then for any $\varepsilon >0$, there exists a positive constant $C_\varepsilon$ such that
\begin{eqnarray}
\label{eq11}
Q(x)|t|^{p-1}|\ln |t|^r|\leq \varepsilon|t|+C_\varepsilon|t|^{l-1},& \forall t\in\mathbb{R}.
\end{eqnarray}
By (\ref{eq11}) and a similar argument in \cite{Squassina 2015}, we obtain that $I\in C^1(\mathcal{H}^{1,2}_0(\Omega),\R)$ and
\begin{eqnarray}
\label{eq12}
    \langle I'(u),v\rangle
&=& a\int_{\Omega\cup \partial \Omega}\nabla u\cdot\nabla vd\mu
      +b\int_{\Omega\cup \partial \Omega} |\nabla u|^2d\mu\int_{\Omega\cup \partial \Omega}\nabla u\cdot\nabla v d\mu\nonumber\\
& & +\lambda \int_\Omega g(x)u^{\frac{2k}{m}-1}v d\mu-\int_\Omega Q(x)|u|^{p-2}u v\ln {|u|^r}d\mu
\end{eqnarray}
for any $u, v\in \mathcal{H}^{1,2}_0(\Omega)$.
\vskip2mm
\noindent
{\bf Definition 2.1.}  {\it $u$ is called as a weak solution of equation (\ref{eq1}) if the following equality holds:
\begin{eqnarray}
\label{e12}
  a\int_{\Omega\cup \partial \Omega}\nabla u\cdot\nabla vd\mu
 +b\int_{\Omega\cup \partial \Omega} |\nabla u|^2d\mu\int_{\Omega\cup \partial \Omega} \nabla u\cdot\nabla vd\mu
 +\int_\Omega \lambda g(x)u^{\frac{2k}{m}-1}v d\mu
=\int_\Omega Q(x)|u|^{p-2}u v\ln {|u|^r} d\mu
\end{eqnarray}
for any $u, v\in \mathcal{H}^{1,2}_0(\Omega)$.
}
\par
By the definition of weak solution, it is easy to see that $u \in \mathcal{H}^{1,2}_0(\Omega)$ is a critical point of $I$ if and only if $u$ is a weak solution of equation (\ref{eq1}).

\vskip2mm
{\section{Proofs}}
\setcounter{equation}{0}
\vskip2mm
\noindent
{\bf Proposition 3.1.}  {\it
If $u\in \mathcal{H}^{1,2}_0(\Omega)$ is a weak solution of equation (\ref{eq1}), then $u$ is  a point-wise solution of (\ref{eq1}).}
\vskip2mm
\noindent
{\bf Proof.} For any fixed $x_0\in \Omega^\circ$, we take a test function $v:\Omega\to \R$ in (\ref{e12}) with
\begin{eqnarray*}
v(x)=\begin{cases}
1,& x=x_0,\\
0,& x\not=x_0.
\end{cases}
\end{eqnarray*}
Then $v\in \mathcal{H}^{1,2}_0(\Omega)$ and
\begin{eqnarray*}
-\left(a+b\int_{\Omega\cup \partial \Omega} |\nabla u(x_0)|^2d\mu\right)\Delta u(x_0)+\lambda g(x_0)u^{\frac{2k}{m}-1}(x_0)-Q(x_0)|u(x_0)|^{p-2}u(x_0)\ln |u(x_0)|^r=0.
\end{eqnarray*}
By the arbitrariness of $x_0$, we conclude that $u$ is  a point-wise solution of (\ref{eq1}).\qed

\vskip2mm
\noindent
{\bf Proposition 3.2.}  {\it  For all $u\in \mathcal{H}^{1,2}_0(\Omega)$, the following equalities hold:}
\begin{eqnarray*}
    I(u)
&=& I(u^+)+I(u^-)-\frac{a}{2}\sum_{x\in\Omega^-}\sum\limits_{y\thicksim x,y\in\Omega^+}w_{xy}u^-(x)u^+(y)
      -\frac{a}{2}\sum_{x\in\Omega^+}\sum\limits_{y\thicksim x,y\in\Omega^-}w_{xy}u^-(y)u^+(x)\\
& & +\frac{b}{4}\left(\sum_{x\in\Omega^-}\sum\limits_{y\thicksim x,y\in\Omega^+}w_{xy}u^-(x)u^+(y)\right)^2
      +\frac{b}{4}\left(\sum_{x\in\Omega^+}\sum\limits_{y\thicksim x,y\in\Omega^-}w_{xy}u^-(y)u^+(x)\right)^2\\
& & +\frac{b}{2}\|u^+\|^2\|u^-\|^2+\frac{b}{2}\sum_{x\in\Omega^-}\sum\limits_{y\thicksim x,y\in\Omega^+}w_{xy}u^-(x)u^+(y)
           \sum_{x\in\Omega^+}\sum\limits_{y\thicksim x,y\in\Omega^-}w_{xy}u^-(y)u^+(x)\\
& &  -\frac{b}{2}\|u^+\|^2\sum_{x\in\Omega^-}\sum\limits_{y\thicksim x,y\in\Omega^+}w_{xy}u^-(x)u^+(y)
       -\frac{b}{2}\|u^+\|^2\sum_{x\in\Omega^+}\sum\limits_{y\thicksim x,y\in\Omega^-}w_{xy}u^-(y)u^+(x)\\
& &  -\frac{b}{2}\|u^-\|^2\sum_{x\in\Omega^-}\sum\limits_{y\thicksim x,y\in\Omega^+}w_{xy}u^-(x)u^+(y)
       -\frac{b}{2}\|u^-\|^2\sum_{x\in\Omega^+}\sum\limits_{y\thicksim x,y\in\Omega^-}w_{xy}u^-(y)u^+(x)
\end{eqnarray*}
and
\begin{eqnarray*}
    \langle I'(u),u^\pm\rangle
&=& \langle I'(u^\pm),u^\pm\rangle-\frac{a}{2}\sum_{x\in\Omega^-}\sum\limits_{y\thicksim x,y\in\Omega^+}w_{xy}u^-(x)u^+(y)
        -\frac{a}{2}\sum_{x\in\Omega^+}\sum\limits_{y\thicksim x,y\in\Omega^-}w_{xy}u^-(y)u^+(x)\\
& & +\frac{b}{2}\left(\sum_{x\in\Omega^-}\sum\limits_{y\thicksim x,y\in\Omega^+}w_{xy}u^-(x)u^+(y)\right)^2
      +\frac{b}{2}\left(\sum_{x\in\Omega^+}\sum\limits_{y\thicksim x,y\in\Omega^-}w_{xy}u^-(y)u^+(x)\right)^2\\
& & +b\|u^+\|^2\|u^-\|^2+b\sum_{x\in\Omega^+}\sum\limits_{y\thicksim x,y\in\Omega^-}w_{xy}u^-(y)u^+(x)
           \sum_{x\in\Omega^-}\sum\limits_{y\thicksim x,y\in\Omega^+}w_{xy}u^-(x)u^+(y)\\
& &  -\frac{3b}{2}\|u^\pm\|^2\sum_{x\in\Omega^-}\sum\limits_{y\thicksim x,y\in\Omega^+}w_{xy}u^-(x)u^+(y)
        -\frac{3b}{2}\|u^\pm\|^2\sum_{x\in\Omega^+}\sum\limits_{y\thicksim x,y\in\Omega^-}w_{xy}u^-(y)u^+(x)\\
& &  -\frac{b}{2}\|u^\mp\|^2\sum_{x\in\Omega^-}\sum\limits_{y\thicksim x,y\in\Omega^+}w_{xy}u^-(x)u^+(y)
        -\frac{b}{2}\|u^\mp\|^2\sum_{x\in\Omega^+}\sum\limits_{y\thicksim x,y\in\Omega^-}w_{xy}u^-(y)u^+(x).
\end{eqnarray*}
{\bf Proof.} Let
\begin{eqnarray*}
\Omega_u^+:=\{x\in\Omega:u(x)>0\},\;\;\Omega_u^-:=\{x\in\Omega:u(x)<0\}\;\;\text{and}\;\;
\mathop{\Omega_u}\limits^\circ:=\{x\in\Omega\cup \partial\Omega:u(x)=0\}.
\end{eqnarray*}
For convenience, we abbreviate $\Omega_u^+$, $\Omega_u^-$ and  $\mathop{\Omega_u}\limits^\circ$ to $\Omega^+$, $\Omega^-$ and  $\mathop{\Omega}\limits^\circ$, respectively. Then, according to the definition of $\Gamma(u,v)$, we have
\begin{eqnarray*}
         \int_{\Omega\cup\partial \Omega} \Gamma(u^+,u^+) d\mu
&  =  &  \int_{\Omega\cup\partial \Omega} |\nabla u^+|^2 d\mu\\
&  =  &  \sum_{x\in \Omega\cup\partial \Omega}\mu(x)\cdot\frac{1}{2\mu(x)}\sum\limits_{y\thicksim x}w_{xy}\big(u^+(y)-u^+(x)\big)^2\\
&  =  &  \frac{1}{2}\sum_{x\in \Omega\cup\partial \Omega}\sum\limits_{y\thicksim x}w_{xy}\big(u^+(y)-u^+(x)\big)^2\\
&  =  &  \frac{1}{2}\sum_{x\in \Omega^+}\sum\limits_{y\thicksim x,y\in\Omega^+}w_{xy}\big(u^+(y)-u^+(x)\big)^2
         +\frac{1}{2}\sum_{x\in \Omega^-}\sum\limits_{y\thicksim x,y\in\Omega^+}w_{xy}\big(u^+(y)-u^+(x)\big)^2\\
&     &  +\frac{1}{2}\sum_{x\in\mathop{\Omega}\limits^\circ}\sum\limits_{y\thicksim x,y\in\Omega^+}w_{xy}\big(u^+(y)-u^+(x)\big)^2
         +\frac{1}{2}\sum_{x\in\Omega^+}\sum\limits_{y\thicksim x,y\in\Omega^-}w_{xy}\big(u^+(y)-u^+(x)\big)^2\\
&     &  +\frac{1}{2}\sum_{x\in\Omega^-}\sum\limits_{y\thicksim x,y\in\Omega^-}w_{xy}\big(u^+(y)-u^+(x)\big)^2
          +\frac{1}{2}\sum_{x\in\mathop{\Omega}\limits^\circ}\sum\limits_{y\thicksim x,y\in\Omega^-}w_{xy}\big(u^+(y)-u^+(x)\big)^2\\
&     &  +\frac{1}{2}\sum_{x\in\Omega^+}\sum\limits_{y\thicksim x,y\in\mathop{\Omega}\limits^\circ}w_{xy}\big(u^+(y)-u^+(x)\big)^2
         +\frac{1}{2}\sum_{x\in\Omega^-}\sum\limits_{y\thicksim x,y\in\mathop{\Omega}\limits^\circ}w_{xy}\big(u^+(y)-u^+(x)\big)^2\\
&     & +\frac{1}{2}\sum_{x\in\mathop{\Omega}\limits^\circ}\sum\limits_{y\thicksim x,y\in\mathop{\Omega}\limits^\circ}w_{xy}\big(u^+(y)-u^+(x)\big)^2\\
&  =  & \frac{1}{2}\sum_{x\in\Omega^+}\sum\limits_{y\thicksim x,y\in\Omega^+}w_{xy}\big(u^+(y)-u^+(x)\big)^2
        +\frac{1}{2}\sum_{x\in\Omega^-}\sum\limits_{y\thicksim x,y\in\Omega^+}w_{xy}\big(u^+(y)\big)^2\\
&     & +\frac{1}{2}\sum_{x\in\mathop{\Omega}\limits^\circ}\sum\limits_{y\thicksim x,y\in\Omega^+}w_{xy}\big(u^+(y)\big)^2
         +\frac{1}{2}\sum_{x\in\Omega^+}\sum\limits_{y\thicksim x,y\in\Omega^-}w_{xy}\big(u^+(x)\big)^2\\
&     & +\frac{1}{2}\sum_{x\in\Omega^+}\sum\limits_{y\thicksim x,y\in\mathop{\Omega}\limits^\circ}w_{xy}\big(u^+(x)\big)^2,
\end{eqnarray*}
and
\begin{eqnarray*}
         \int_{\Omega\cup\partial \Omega} \Gamma(u^-,u^-) d\mu
&  =  &  \int_{\Omega\cup\partial \Omega} |\nabla u^-|^2 d\mu\\
&  =  &  \frac{1}{2}\sum_{x\in \Omega\cup\partial \Omega}\sum\limits_{y\thicksim x}w_{xy}\big(u^-(y)-u^-(x)\big)^2\\
&  =  &  \frac{1}{2}\sum_{x\in \Omega^-}\sum\limits_{y\thicksim x,y\in\Omega^+}w_{xy}\big(u^-(x)\big)^2
         +\frac{1}{2}\sum_{x\in\Omega^+}\sum\limits_{y\thicksim x,y\in\Omega^-}w_{xy}\big(u^-(y)\big)^2\\
&     &  +\frac{1}{2}\sum_{x\in\Omega^-}\sum\limits_{y\thicksim x,y\in\Omega^-}w_{xy}\big(u^-(y)-u^-(x)\big)^2
          +\frac{1}{2}\sum_{x\in\mathop{\Omega}\limits^\circ}\sum\limits_{y\thicksim x,y\in\Omega^-}w_{xy}\big(u^-(y)\big)^2\\
&     &  +\frac{1}{2}\sum_{x\in\Omega^-}\sum\limits_{y\thicksim x,y\in\mathop{\Omega}\limits^\circ}w_{xy}\big(u^-(x)\big)^2.
\end{eqnarray*}
Moreover,
\begin{eqnarray*}
          \int_{\Omega\cup\partial \Omega}|\nabla u|^2d\mu
&  =  &  \frac{1}{2}\sum_{x\in\Omega\cup\partial \Omega}\sum\limits_{y\thicksim x}w_{xy}\big(u(y)-u(x)\big)^2\\
&  =  &  \frac{1}{2}\sum_{x\in \Omega^+}\sum\limits_{y\thicksim x,y\in\Omega^+}w_{xy}\big(u^+(y)-u^+(x)\big)^2
         +\frac{1}{2}\sum_{x\in \Omega^-}\sum\limits_{y\thicksim x,y\in\Omega^+}w_{xy}\big(u^+(y)-u^-(x)\big)^2\\
&     &  +\frac{1}{2}\sum_{x\in\mathop{\Omega}\limits^\circ}\sum\limits_{y\thicksim x,y\in\Omega^+}w_{xy}\big(u^+(y)\big)^2
         +\frac{1}{2}\sum_{x\in\Omega^+}\sum\limits_{y\thicksim x,y\in\Omega^-}w_{xy}\big(u^-(y)-u^+(x)\big)^2\\
&     &  +\frac{1}{2}\sum_{x\in\Omega^-}\sum\limits_{y\thicksim x,y\in\Omega^-}w_{xy}\big(u^-(y)-u^-(x)\big)^2
          +\frac{1}{2}\sum_{x\in\mathop{\Omega}\limits^\circ}\sum\limits_{y\thicksim x,y\in\Omega^-}w_{xy}\big(u^-(y)\big)^2\\
&     &  +\frac{1}{2}\sum_{x\in\Omega^+}\sum\limits_{y\thicksim x,y\in\mathop{\Omega}\limits^\circ}w_{xy}\big(u^+(x)\big)^2
         +\frac{1}{2}\sum_{x\in\Omega^-}\sum\limits_{y\thicksim x,y\in\mathop{\Omega}\limits^\circ}w_{xy}\big(u^-(x)\big)^2\\
&  =  &  \frac{1}{2}\sum_{x\in \Omega^+}\sum\limits_{y\thicksim x,y\in\Omega^+}w_{xy}\big(u^+(y)-u^+(x)\big)^2
         +\frac{1}{2}\sum_{x\in \Omega^-}\sum\limits_{y\thicksim x,y\in\Omega^+}w_{xy}\big(u^+(y)\big)^2\\
&     &  -\sum_{x\in \Omega^-}\sum\limits_{y\thicksim x,y\in\Omega^+}w_{xy}u^+(y)u^-(x)
         +\frac{1}{2}\sum_{x\in \Omega^-}\sum\limits_{y\thicksim x,y\in\Omega^+}w_{xy}\big(u^-(x)\big)^2\\
&     &  +\frac{1}{2}\sum_{x\in\mathop{\Omega}\limits^\circ}\sum\limits_{y\thicksim x,y\in\Omega^+}w_{xy}\big(u^+(y)\big)^2
         +\frac{1}{2}\sum_{x\in\Omega^+}\sum\limits_{y\thicksim x,y\in\Omega^-}w_{xy}\big(u^-(y)\big)^2\\
&     &  -\sum_{x\in\Omega^+}\sum\limits_{y\thicksim x,y\in\Omega^-}w_{xy}u^-(y)u^+(x)
         +\frac{1}{2}\sum_{x\in\Omega^+}\sum\limits_{y\thicksim x,y\in\Omega^-}w_{xy}\big(u^+(x)\big)^2\\
&     &  +\frac{1}{2}\sum_{x\in\Omega^-}\sum\limits_{y\thicksim x,y\in\Omega^-}w_{xy}\big(u^-(y)-u^-(x)\big)^2
         +\frac{1}{2}\sum_{x\in\mathop{\Omega}\limits^\circ}\sum\limits_{y\thicksim x,y\in\Omega^-}w_{xy}\big(u^-(y)\big)^2\\
&     &  +\frac{1}{2}\sum_{x\in\Omega^+}\sum\limits_{y\thicksim x,y\in\mathop{\Omega}\limits^\circ}w_{xy}\big(u^+(x)\big)^2
         +\frac{1}{2}\sum_{x\in\Omega^-}\sum\limits_{y\thicksim x,y\in\mathop{\Omega}\limits^\circ}w_{xy}\big(u^-(x)\big)^2.
\end{eqnarray*}
So, we have
\begin{small}
\begin{eqnarray*}
 \int_{\Omega\cup\partial \Omega} |\nabla u|^2 d\mu
=\int_{\Omega\cup\partial \Omega}|\nabla u^+|^2d\mu+\int_{\Omega\cup\partial \Omega}|\nabla u^-|^2 d\mu
      -\sum_{x\in\Omega^-}\sum\limits_{y\thicksim x,y\in\Omega^+}w_{xy}u^-(x)u^+(y)
      -\sum_{x\in\Omega^+}\sum\limits_{y\thicksim x,y\in\Omega^-}w_{xy}u^-(y)u^+(x),
\end{eqnarray*}
\end{small}that is,
\begin{eqnarray}
\label{g1}
 \|u\|^2
=\|u^+\|^2 + \|u^-\|^2-\sum_{x\in\Omega^-}\sum\limits_{y\thicksim x,y\in\Omega^+}w_{xy}u^-(x)u^+(y)
              -\sum_{x\in\Omega^+}\sum\limits_{y\thicksim x,y\in\Omega^-}w_{xy}u^-(y)u^+(x).
\end{eqnarray}
In addition, according to (\ref{g1}), we can obtain that
\begin{eqnarray}
\label{g2}
         \left(\int_{\Omega\cup\partial \Omega} |\nabla u|^2 d\mu\right)^2
&  =  &  \left(\int_{\Omega\cup\partial \Omega}|\nabla u^+|^2d\mu+ \int_{\Omega\cup\partial \Omega}|\nabla u^-|^2 d\mu\right.\nonumber\\
&     &   \left.-\sum_{x\in\Omega^-}\sum\limits_{y\thicksim x,y\in\Omega^+}w_{xy}u^-(x)u^+(y)
                        -\sum_{x\in\Omega^+}\sum\limits_{y\thicksim x,y\in\Omega^-}w_{xy}u^-(y)u^+(x)\right)^2\nonumber\\
&  =  &  \left(\int_{\Omega\cup\partial \Omega}|\nabla u^+|^2d\mu\right)^2
             +\left(\int_{\Omega\cup\partial \Omega}|\nabla u^-|^2 d\mu\right)^2
             +\left(\sum_{x\in\Omega^-}\sum\limits_{y\thicksim x,y\in\Omega^+}w_{xy}u^-(x)u^+(y)\right)^2\nonumber\\
&     &  +\left(\sum_{x\in\Omega^+}\sum\limits_{y\thicksim x,y\in\Omega^-}w_{xy}u^-(y)u^+(x)\right)^2
             +2\int_{\Omega\cup\partial \Omega} |\nabla u^+|^2 d\mu\int_{\Omega\cup\partial \Omega}|\nabla u^-|^2d\mu\nonumber\\
&     &  -2\sum_{x\in\Omega^-}\sum\limits_{y\thicksim x,y\in\Omega^+}w_{xy}u^-(x)u^+(y)\int_{\Omega\cup\partial \Omega} |\nabla u^+|^2 d\mu\nonumber\\
&    &   -2\sum_{x\in\Omega^+}\sum\limits_{y\thicksim x,y\in\Omega^-}w_{xy}u^-(y)u^+(x)\int_{\Omega\cup\partial \Omega} |\nabla u^+|^2 d\mu\nonumber\\
&     &  -2\sum_{x\in\Omega^-}\sum\limits_{y\thicksim x,y\in\Omega^+}w_{xy}u^-(x)u^+(y)\int_{\Omega\cup\partial \Omega} |\nabla u^-|^2 d\mu\nonumber\\
&    &   -2\sum_{x\in\Omega^+}\sum\limits_{y\thicksim x,y\in\Omega^-}w_{xy}u^-(y)u^+(x)\int_{\Omega\cup\partial \Omega} |\nabla u^-|^2 d\mu\nonumber\\
&    &   +2\sum_{x\in\Omega^-}\sum\limits_{y\thicksim x,y\in\Omega^+}w_{xy}u^-(x)u^+(y)\sum_{x\in\Omega^+}\sum\limits_{y\thicksim x,y\in\Omega^-}w_{xy}u^-(y)u^+(x).
\end{eqnarray}
Similarly, we can also obtain that
\begin{eqnarray*}
         \int_{\Omega\cup\partial \Omega}\nabla u\cdot\nabla u^+d\mu
&  =  &\frac{1}{2}\sum_{x\in\Omega\cup\partial \Omega}\sum\limits_{y\thicksim x}w_{xy}(u(y)-u(x))(u^+(y)-u^+(x))\\
&  =  &\frac{1}{2}\sum_{x\in \Omega^+}\sum\limits_{y\thicksim x,y\in\Omega^+}w_{xy}(u^+(y)-u^+(x))^2+\frac{1}{2}\sum_{x\in \Omega^-}\sum\limits_{y\thicksim x,y\in\Omega^+}w_{xy}(u^+(y)-u^-(x))(u^+(y))\\
&     &+\frac{1}{2}\sum_{x\in\mathop{\Omega}\limits^\circ}\sum\limits_{y\thicksim x,y\in\Omega^+}w_{xy}(u^+(y))^2
       +\frac{1}{2}\sum_{x\in \Omega^+}\sum\limits_{y\thicksim x,y\in\Omega^-}w_{xy}(u^-(y)-u^+(x))(-u^+(x))\\
&     &+\frac{1}{2}\sum_{x\in \Omega^+}\sum\limits_{y\thicksim x,y\in\mathop{\Omega}\limits^\circ}w_{xy}(-u^+(x))(-u^+(x))\\
&  =  &\frac{1}{2}\sum_{x\in \Omega^+}\sum\limits_{y\thicksim x,y\in\Omega^+}w_{xy}(u^+(y)-u^+(x))^2
       +\frac{1}{2}\sum_{x\in \Omega^-}\sum\limits_{y\thicksim x,y\in\Omega^+}w_{xy}(u^+(y))^2\\
&     &-\frac{1}{2}\sum_{x\in \Omega^-}\sum\limits_{y\thicksim x,y\in\Omega^+}w_{xy}u^-(x)u^+(y)
       +\frac{1}{2}\sum_{x\in\mathop{\Omega}\limits^\circ}\sum\limits_{y\thicksim x,y\in\Omega^+}w_{xy}(u^+(y))^2\\
&     &-\frac{1}{2}\sum_{x\in \Omega^+}\sum\limits_{y\thicksim x,y\in\Omega^-}w_{xy}u^-(y)u^+(x)
       +\frac{1}{2}\sum_{x\in \Omega^+}\sum\limits_{y\thicksim x,y\in\Omega^-}w_{xy}(u^+(x))^2\\
&     &+\frac{1}{2}\sum_{x\in \Omega^+}\sum\limits_{y\thicksim x,y\in\mathop{\Omega}\limits^\circ}w_{xy}(u^+(x))^2,
\end{eqnarray*}
and
\begin{eqnarray*}
         \int_{\Omega\cup\partial \Omega}\nabla u\cdot\nabla u^-d\mu
&  =  &\frac{1}{2}\sum_{x\in\Omega\cup\partial \Omega}\sum\limits_{y\thicksim x}w_{xy}(u(y)-u(x))(u^-(y)-u^-(x))\\
&  =  &\frac{1}{2}\sum_{x\in \Omega^-}\sum\limits_{y\thicksim x,y\in\Omega^+}w_{xy}(u^+(y)-u^-(x))(-u^-(x))\\
&     &+\frac{1}{2}\sum_{x\in \Omega^+}\sum\limits_{y\thicksim x,y\in\Omega^-}w_{xy}(u^-(y)-u^+(x))(u^-(y))\\
&     &+\frac{1}{2}\sum_{x\in \Omega^-}\sum\limits_{y\thicksim x,y\in\Omega^-}w_{xy}(u^-(y)-u^-(x))^2
       +\frac{1}{2}\sum_{x\in\mathop{\Omega}\limits^\circ}\sum\limits_{y\thicksim x,y\in\Omega^-}w_{xy}(u^-(y))^2\\
&     &+\frac{1}{2}\sum_{x\in \Omega^-}\sum\limits_{y\thicksim x,y\in\mathop{\Omega}\limits^\circ}w_{xy}(-u^-(x))(-u^-(x))\\
&  =  &-\frac{1}{2}\sum_{x\in \Omega^-}\sum\limits_{y\thicksim x,y\in\Omega^+}w_{xy}u^+(y)u^-(x)
       +\frac{1}{2}\sum_{x\in \Omega^-}\sum\limits_{y\thicksim x,y\in\Omega^+}w_{xy}(u^-(x))^2\\
&     &+\frac{1}{2}\sum_{x\in \Omega^+}\sum\limits_{y\thicksim x,y\in\Omega^-}w_{xy}(u^-(y))^2
       -\frac{1}{2}\sum_{x\in \Omega^+}\sum\limits_{y\thicksim x,y\in\Omega^-}w_{xy}u^+(x)u^-(y)\\
&     &+\frac{1}{2}\sum_{x\in \Omega^-}\sum\limits_{y\thicksim x,y\in\Omega^-}w_{xy}(u^-(y)-u^-(x))^2
       +\frac{1}{2}\sum_{x\in\mathop{\Omega}\limits^\circ}\sum\limits_{y\thicksim x,y\in\Omega^-}w_{xy}(u^-(y))^2\\
&     &+\frac{1}{2}\sum_{x\in \Omega^-}\sum\limits_{y\thicksim x,y\in\mathop{\Omega}\limits^\circ}w_{xy}(u^-(x))^2.
\end{eqnarray*}
Then
\begin{small}
\begin{eqnarray}
\label{g3}
   \int_{\Omega\cup\partial \Omega}\nabla u\cdot\nabla u^+d\mu
=  \int_{\Omega\cup\partial \Omega}|\nabla u^+|^2d\mu
     -\frac{1}{2}\sum_{x\in\Omega^-}\sum\limits_{y\thicksim x,y\in\Omega^+}w_{xy}u^-(x)u^+(y)
     -\frac{1}{2}\sum_{x\in\Omega^+}\sum\limits_{y\thicksim x,y\in\Omega^-}w_{xy}u^-(y)u^+(x)
\end{eqnarray}
\end{small}and
\begin{small}
\begin{eqnarray}
\label{g4}
   \int_{\Omega\cup\partial \Omega}\nabla u\cdot\nabla u^-d\mu
=  \int_{\Omega\cup\partial \Omega}|\nabla u^-|^2d\mu
     -\frac{1}{2}\sum_{x\in\Omega^-}\sum\limits_{y\thicksim x,y\in\Omega^+}w_{xy}u^-(x)u^+(y)
     -\frac{1}{2}\sum_{x\in\Omega^+}\sum\limits_{y\thicksim x,y\in\Omega^-}w_{xy}u^-(y)u^+(x).
\end{eqnarray}
\end{small}In virtue of (\ref{g1}), (\ref{g3}) and (\ref{g4}), we have
\begin{eqnarray}
\label{g5}
& & \int_{\Omega\cup \partial \Omega} |\nabla u|^2d\mu\int_{\Omega\cup \partial \Omega}\nabla u\cdot\nabla u^\pm d\mu\nonumber\\
&=& \left(\int_{\Omega\cup\partial \Omega}|\nabla u^+|^2d\mu+ \int_{\Omega\cup\partial \Omega}|\nabla u^-|^2 d\mu
       -\sum_{x\in\Omega^-}\sum\limits_{y\thicksim x,y\in\Omega^+}w_{xy}u^-(x)u^+(y)
       -\sum_{x\in\Omega^+}\sum\limits_{y\thicksim x,y\in\Omega^-}w_{xy}u^-(y)u^+(x)\right)\nonumber\\
& &  \times\left(\int_{\Omega\cup\partial \Omega}|\nabla u^\pm|^2d\mu
        -\frac{1}{2}\sum_{x\in\Omega^-}\sum\limits_{y\thicksim x,y\in\Omega^+}w_{xy}u^-(x)u^+(y)
        -\frac{1}{2}\sum_{x\in\Omega^+}\sum\limits_{y\thicksim x,y\in\Omega^-}w_{xy}u^-(y)u^+(x)\right)\nonumber\\
&=&  \int_{\Omega\cup\partial \Omega}|\nabla u^+|^2d\mu\int_{\Omega\cup\partial \Omega}|\nabla u^\pm|^2d\mu
        +\int_{\Omega\cup\partial \Omega}|\nabla u^-|^2d\mu\int_{\Omega\cup\partial \Omega}|\nabla u^\pm|^2d\mu\nonumber\\
& &  -\sum_{x\in\Omega^-}\sum\limits_{y\thicksim x,y\in\Omega^+}w_{xy}u^-(x)u^+(y)
           \int_{\Omega\cup\partial \Omega}|\nabla u^\pm|^2d\mu
        -\sum_{x\in\Omega^+}\sum\limits_{y\thicksim x,y\in\Omega^-}w_{xy}u^-(y)u^+(x)
           \int_{\Omega\cup\partial \Omega}|\nabla u^\pm|^2d\mu\nonumber\\
& &  -\frac{1}{2}\sum_{x\in\Omega^-}\sum\limits_{y\thicksim x,y\in\Omega^+}w_{xy}u^-(x)u^+(y)
           \int_{\Omega\cup\partial \Omega}|\nabla u^+|^2d\mu
        -\frac{1}{2}\sum_{x\in\Omega^-}\sum\limits_{y\thicksim x,y\in\Omega^+}w_{xy}u^-(x)u^+(y)
           \int_{\Omega\cup\partial \Omega}|\nabla u^-|^2d\mu\nonumber\\
& &  +\frac{1}{2}\left(\sum_{x\in\Omega^-}\sum\limits_{y\thicksim x,y\in\Omega^+}w_{xy}u^-(x)u^+(y)\right)^2
        +\frac{1}{2}\sum_{x\in\Omega^-}\sum\limits_{y\thicksim x,y\in\Omega^+}w_{xy}u^-(x)u^+(y)
           \sum_{x\in\Omega^+}\sum\limits_{y\thicksim x,y\in\Omega^-}w_{xy}u^-(y)u^+(x)\nonumber\\
& &  -\frac{1}{2}\sum_{x\in\Omega^+}\sum\limits_{y\thicksim x,y\in\Omega^-}w_{xy}u^-(y)u^+(x)
           \int_{\Omega\cup\partial \Omega}|\nabla u^+|^2d\mu
        -\frac{1}{2}\sum_{x\in\Omega^+}\sum\limits_{y\thicksim x,y\in\Omega^-}w_{xy}u^-(y)u^+(x)
           \int_{\Omega\cup\partial \Omega}|\nabla u^-|^2d\mu\nonumber\\
& &  +\frac{1}{2}\sum_{x\in\Omega^+}\sum\limits_{y\thicksim x,y\in\Omega^-}w_{xy}u^-(y)u^+(x)
           \sum_{x\in\Omega^-}\sum\limits_{y\thicksim x,y\in\Omega^+}w_{xy}u^-(x)u^+(y)
         +\frac{1}{2}\left(\sum_{x\in\Omega^+}\sum\limits_{y\thicksim x,y\in\Omega^-}w_{xy}u^-(y)u^+(x)\right)^2\nonumber\\
&=&  \left(\int_{\Omega\cup\partial \Omega}|\nabla u^\pm|^2d\mu\right)^2
         +\int_{\Omega\cup\partial \Omega}|\nabla u^+|^2d\mu\int_{\Omega\cup\partial \Omega}|\nabla u^-|^2d\mu           +\frac{1}{2}\left(\sum_{x\in\Omega^-}\sum\limits_{y\thicksim x,y\in\Omega^+}w_{xy}u^-(x)u^+(y)\right)^2\nonumber\\
& &  +\frac{1}{2}\left(\sum_{x\in\Omega^+}\sum\limits_{y\thicksim x,y\in\Omega^-}w_{xy}u^-(y)u^+(x)\right)^2
         +\sum_{x\in\Omega^+}\sum\limits_{y\thicksim x,y\in\Omega^-}w_{xy}u^-(y)u^+(x)
           \sum_{x\in\Omega^-}\sum\limits_{y\thicksim x,y\in\Omega^+}w_{xy}u^-(x)u^+(y)\nonumber\\
& &  -\frac{3}{2}\sum_{x\in\Omega^-}\sum\limits_{y\thicksim x,y\in\Omega^+}w_{xy}u^-(x)u^+(y)
           \int_{\Omega\cup\partial \Omega}|\nabla u^\pm|^2d\mu
        -\frac{3}{2}\sum_{x\in\Omega^+}\sum\limits_{y\thicksim x,y\in\Omega^-}w_{xy}u^-(y)u^+(x)
           \int_{\Omega\cup\partial \Omega}|\nabla u^\pm|^2d\mu\nonumber\\
& &  -\frac{1}{2}\sum_{x\in\Omega^-}\sum\limits_{y\thicksim x,y\in\Omega^+}w_{xy}u^-(x)u^+(y)
           \int_{\Omega\cup\partial \Omega}|\nabla u^\mp|^2d\mu
        -\frac{1}{2}\sum_{x\in\Omega^+}\sum\limits_{y\thicksim x,y\in\Omega^-}w_{xy}u^-(y)u^+(x)
           \int_{\Omega\cup\partial \Omega}|\nabla u^\mp|^2d\mu.\nonumber\\
\end{eqnarray}
By (\ref{g1}), (\ref{g2}) and (\ref{g5}), it is easy to see that the conclusions hold.
  \qed

\vskip2mm
\par
Next, we establish an  inequality related to $I(u)$,\;\;$I(su^++tu^-)$,\;\;$\langle I'(u),u^+\rangle$
and $\langle I'(u),u^-\rangle$.
\vskip2mm
\noindent
{\bf Lemma 3.1.}  {\it For all $u\in \mathcal{H}^{1,2}_0(\Omega)$ and $s,t\geq0$, there exists}
\begin{eqnarray}
\label{Eq4}
      I(u)
&\geq& I(su^++tu^-)+\frac{1-s^{p}}{p}\langle I'(u),u^+\rangle
         +\frac{1-t^{p}}{p}\langle I'(u),u^-\rangle+a\left(\frac{1-s^2}{2}-\frac{1-s^{p}}{p}\right)\|u^+\|^2\nonumber\\
&    & +a\left(\frac{1-t^2}{2}-\frac{1-t^{p}}{p}\right)\|u^-\|^2
        +b\left(\frac{1-s^4}{4}-\frac{1-s^{p}}{p}\right)\|u^+\|^4
        +b\left(\frac{1-t^4}{4}-\frac{1-t^{p}}{p}\right)\|u^-\|^4\nonumber\\
&    & +\frac{b(s^2-t^2)^2}{4}\|u^+\|^2\|u^-\|^2
        +\frac{b(s^2-t^2)^2}{4}\sum_{x\in\Omega^-}\sum\limits_{y\thicksim x,y\in\Omega^+}w_{xy}u^-(x)u^+(y)
           \sum_{x\in\Omega^+}\sum\limits_{y\thicksim x,y\in\Omega^-}w_{xy}u^-(y)u^+(x)\nonumber\\
&    &-\frac{a(s-t)^2}{4}\sum_{x\in\Omega^-}\sum\limits_{y\thicksim x,y\in\Omega^+}w_{xy}u^-(x)u^+(y)
          -\frac{a(s-t)^2}{4}\sum_{x\in\Omega^+}\sum\limits_{y\thicksim x,y\in\Omega^-}w_{xy}u^-(y)u^+(x)\nonumber\\
&    &+\frac{b(s^2-t^2)^2}{8}\left(\sum_{x\in\Omega^-}\sum\limits_{y\thicksim x,y\in\Omega^+}w_{xy}u^-(x)u^+(y)\right)^2
+\frac{b(s^2-t^2)^2}{8}\left(\sum_{x\in\Omega^+}\sum\limits_{y\thicksim x,y\in\Omega^-}w_{xy}u^-(y)u^+(x)\right)^2.\nonumber\\
\end{eqnarray}
{\bf Proof.} It is easy to see that (\ref{Eq4}) holds for $u=0$. Next, we let $u\neq 0$. By Appendix 1 below, there holds
\begin{eqnarray}
\label{Eq5}
r(1-\tau^{p})+p\tau^{p}\ln\tau^r> 0,\;\;\forall \tau\in(0,1)\cup(1,+\infty).
\end{eqnarray}
For $u\in \mathcal{H}^{1,2}_0(\Omega)\backslash\{0\}$ and all $s,t\ge 0$, we have
\begin{eqnarray}
\label{Eq6}
&&\int_\Omega |su^++tu^-|^{p}\ln |su^++tu^-|^r d\mu\nonumber\\
&=&\int_{\Omega^+} |su^++tu^-|^{p}\ln |su^++tu^-|^r d\mu+\int_{\Omega^-} |su^++tu^-|^{p}\ln |su^++tu^-|^r d\mu\nonumber\\
&=&\int_{\Omega^+} |su^+|^{p}\ln |su^+|^r d\mu+\int_{\Omega^-} |tu^-|^{p}\ln |tu^-|^r d\mu\nonumber\\
&=&\int_\Omega \Big[|su^+|^{p}\left(\ln |u^+|^r+\ln s^r\right)+|tu^-|^{p}\left(\ln |u^-|^r+\ln t^r\right)\Big] d\mu.
\end{eqnarray}
Because the function $f(x)=\frac{1-a^x}{x}$ is monotonically decreasing on $(0,+\infty)$ for $a>0$ and $a\neq 1$ \big(see Appendix 2\big), in combination with (\ref{eq9}), (\ref{eq12}), (\ref{g1}), (\ref{g2}), (\ref{Eq5}) and (\ref{Eq6}), we have
\begin{eqnarray}
\label{Eq7}
&    &I(u)-I(su^++tu^-)\nonumber\\
& =  &\frac{a}{2}\left(\|u\|^2-\|su^++tu^-\|^2\right)+\frac{b}{4}\left(\|u\|^4-\|su^++tu^-\|^4\right)
      +\frac{m}{2k}\int_\Omega\lambda g(x)\left[u^{\frac{2k}{m}}-(su^++tu^-)^{\frac{2k}{m}}\right]d\mu\nonumber\\
&    &+\frac{r}{p^2}\int_\Omega Q(x)\left(|u|^p-|su^++tu^-|^p\right)d\mu
      -\frac{1}{p}\int_\Omega Q(x)\Big(|u|^{p}\ln |u|^r-|su^++tu^-|^{p}\ln |su^++tu^-|^r\Big)d\mu\nonumber\\
& =  &\frac{a(1-s^2)}{2}\|u^+\|^2+\frac{a(1-t^2)}{2}\|u^-\|^2+\frac{b(1-s^4)}{4}\|u^+\|^4+\frac{b(1-t^4)}{4}\|u^-\|^4\nonumber\\
&    &+\frac{a(st-1)}{2}\sum_{x\in\Omega^-}\sum\limits_{y\thicksim x,y\in\Omega^+}w_{xy}u^-(x)u^+(y)
      +\frac{a(st-1)}{2}\sum_{x\in\Omega^+}\sum\limits_{y\thicksim x,y\in\Omega^-}w_{xy}u^-(y)u^+(x)\nonumber\\
&    &+\frac{b(1-s^2t^2)}{4}\left(\sum_{x\in\Omega^-}\sum\limits_{y\thicksim x,y\in\Omega^+}w_{xy}u^-(x)u^+(y)\right)^2
      +\frac{b(1-s^2t^2)}{4}\left(\sum_{x\in\Omega^+}\sum\limits_{y\thicksim x,y\in\Omega^-}w_{xy}u^-(y)u^+(x)\right)^2\nonumber\\
&    &+\frac{b(1-s^2t^2)}{2}\|u^+\|^2\|u^-\|^2
      +\frac{b(1-s^2t^2)}{2}\sum_{x\in\Omega^-}\sum\limits_{y\thicksim x,y\in\Omega^+}w_{xy}u^-(x)u^+(y)
           \sum_{x\in\Omega^+}\sum\limits_{y\thicksim x,y\in\Omega^-}w_{xy}u^-(y)u^+(x)\nonumber\\
&    &+\frac{b(s^3t-1)}{2}\|u^+\|^2\sum_{x\in\Omega^-}\sum\limits_{y\thicksim x,y\in\Omega^+}w_{xy}u^-(x)u^+(y)
      +\frac{b(s^3t-1)}{2}\|u^+\|^2\sum_{x\in\Omega^+}\sum\limits_{y\thicksim x,y\in\Omega^-}w_{xy}u^-(y)u^+(x)\nonumber\\
&    &+\frac{b(st^3-1)}{2}\|u^-\|^2\sum_{x\in\Omega^-}\sum\limits_{y\thicksim x,y\in\Omega^+}w_{xy}u^-(x)u^+(y)
      +\frac{b(st^3-1)}{2}\|u^-\|^2\sum_{x\in\Omega^+}\sum\limits_{y\thicksim x,y\in\Omega^-}w_{xy}u^-(y)u^+(x)\nonumber\\
&    &+\frac{r(1-s^p)}{p^2}\int_\Omega Q(x)|u^+|^pd\mu+\frac{r(1-t^p)}{p^2}\int_\Omega Q(x)|u^-|^pd\mu\nonumber\\
&    &-\frac{1-s^{p}}{p}\int_\Omega Q(x)|u^+|^{p}\ln|u^+|^rd\mu
      -\frac{1-t^{p}}{p}\int_\Omega Q(x)|u^-|^{p}\ln|u^-|^rd\mu
      +\frac{\lambda(1-s^{\frac{2k}{m}})}{\frac{2k}{m}}\int_\Omega g(x)(u^+)^{\frac{2k}{m}}d\mu\nonumber\\
&    &+\frac{\lambda(1-t^{\frac{2k}{m}})}{\frac{2k}{m}}\int_\Omega g(x)(u^-)^{\frac{2k}{m}}d\mu
      +\frac{1}{p}\int_\Omega Q(x)|su^+|^{p}\ln s^rd\mu+\frac{1}{p}\int_\Omega Q(x)|tu^-|^{p}\ln t^rd\mu\nonumber\\
& =  &\frac{1-s^{p}}{p}\left[\langle I'(u),u^+\rangle-a\int_{\Omega\cup\partial\Omega}\nabla u\cdot\nabla u^+d\mu
        -b\int_{\Omega\cup\partial\Omega}|\nabla u|^2d\mu\int_{\Omega\cup\partial\Omega}\nabla u\cdot\nabla u^+d\mu
        -\lambda\int_\Omega g(x)(u^+)^{\frac{2k}{m}}d\mu\right]\nonumber\\
&    &+\frac{1-t^{p}}{p}\left[\langle I'(u),u^-\rangle-a\int_{\Omega\cup\partial\Omega}\nabla u\cdot\nabla u^-d\mu
        -b\int_{\Omega\cup\partial\Omega}|\nabla u|^2d\mu\int_{\Omega\cup\partial\Omega}\nabla u\cdot\nabla u^-d\mu
        -\lambda\int_\Omega g(x)(u^-)^{\frac{2k}{m}}d\mu\right]\nonumber\\
&    &+\frac{a(1-s^2)}{2}\|u^+\|^2+\frac{a(1-t^2)}{2}\|u^-\|^2
                 +\frac{b(1-s^4)}{4}\|u^+\|^4+\frac{b(1-t^4)}{4}\|u^-\|^4\nonumber\\
&    &+\frac{a(st-1)}{2}\sum_{x\in\Omega^-}\sum\limits_{y\thicksim x,y\in\Omega^+}w_{xy}u^-(x)u^+(y)
      +\frac{a(st-1)}{2}\sum_{x\in\Omega^+}\sum\limits_{y\thicksim x,y\in\Omega^-}w_{xy}u^-(y)u^+(x)\nonumber\\
&    &+\frac{b(1-s^2t^2)}{4}\left(\sum_{x\in\Omega^-}\sum\limits_{y\thicksim x,y\in\Omega^+}w_{xy}u^-(x)u^+(y)\right)^2
      +\frac{b(1-s^2t^2)}{4}\left(\sum_{x\in\Omega^+}\sum\limits_{y\thicksim x,y\in\Omega^-}w_{xy}u^-(y)u^+(x)\right)^2\nonumber\\
&    &+\frac{b(1-s^2t^2)}{2}\|u^+\|^2\|u^-\|^2
      +\frac{b(1-s^2t^2)}{2}\sum_{x\in\Omega^-}\sum\limits_{y\thicksim x,y\in\Omega^+}w_{xy}u^-(x)u^+(y)
           \sum_{x\in\Omega^+}\sum\limits_{y\thicksim x,y\in\Omega^-}w_{xy}u^-(y)u^+(x)\nonumber\\
&    &+\frac{b(s^3t-1)}{2}\|u^+\|^2\sum_{x\in\Omega^-}\sum\limits_{y\thicksim x,y\in\Omega^+}w_{xy}u^-(x)u^+(y)
      +\frac{b(s^3t-1)}{2}\|u^+\|^2\sum_{x\in\Omega^+}\sum\limits_{y\thicksim x,y\in\Omega^-}w_{xy}u^-(y)u^+(x)\nonumber\\
&    &+\frac{b(st^3-1)}{2}\|u^-\|^2\sum_{x\in\Omega^-}\sum\limits_{y\thicksim x,y\in\Omega^+}w_{xy}u^-(x)u^+(y)
      +\frac{b(st^3-1)}{2}\|u^-\|^2\sum_{x\in\Omega^+}\sum\limits_{y\thicksim x,y\in\Omega^-}w_{xy}u^-(y)u^+(x)\nonumber\\
&    &+\frac{r(1-s^p)}{p^2}\int_\Omega Q(x)|u^+|^pd\mu+\frac{r(1-t^p)}{p^2}\int_\Omega Q(x)|u^-|^pd\mu
      +\frac{\lambda(1-s^{\frac{2k}{m}})}{\frac{2k}{m}}\int_\Omega g(x)(u^+)^{\frac{2k}{m}}d\mu\nonumber\\
&    &+\frac{\lambda(1-t^{\frac{2k}{m}})}{\frac{2k}{m}}\int_\Omega g(x)(u^-)^{\frac{2k}{m}}d\mu
      +\frac{1}{p}\int_\Omega Q(x)|su^+|^{p}\ln s^rd\mu+\frac{1}{p}\int_\Omega Q(x)|tu^-|^{p}\ln t^rd\mu\nonumber\\
& =  &\frac{1-s^{p}}{p}\langle I'(u),u^+\rangle+\frac{1-t^{p}}{p}\langle I'(u),u^-\rangle
        +a\left(\frac{1-s^2}{2}-\frac{1-s^{p}}{p}\right)\|u^+\|^2
        +a\left(\frac{1-t^2}{2}-\frac{1-t^{p}}{p}\right)\|u^-\|^2\nonumber\\
&    &+b\left(\frac{1-s^4}{4}-\frac{1-s^{p}}{p}\right)\|u^+\|^4
        +b\left(\frac{1-t^4}{4}-\frac{1-t^{p}}{p}\right)\|u^-\|^4
        +\lambda\left(\frac{1-s^{\frac{2k}{m}}}{\frac{2k}{m}}-\frac{1-s^{p}}{p}\right)\int_\Omega g(x)(u^+)^{\frac{2k}{m}}d\mu\nonumber\\
&    &+\lambda\left(\frac{1-t^{\frac{2k}{m}}}{\frac{2k}{m}}-\frac{1-t^{p}}{p}\right)\int_\Omega g(x)(u^-)^{\frac{2k}{m}}d\mu
      +\frac{r(1-s^{p})+ps^{p}\ln s^r}{p^2}\int_\Omega Q(x)|u^+|^{p}d\mu\nonumber\\
&    &+\frac{r(1-t^{p})+pt^{p}\ln t^r}{p^2}\int_\Omega Q(x)|u^-|^{p}d\mu
        +b\left(\frac{1-s^2t^2}{2}-\frac{1-s^{p}}{p}-\frac{1-t^{p}}{p}\right)\|u^+\|^2\|u^-\|^2\nonumber\\
&    &+a\left(\frac{st-1}{2}+\frac{1-s^{p}}{2p}+\frac{1-t^{p}}{2p}\right)\sum_{x\in\Omega^-}
          \sum\limits_{y\thicksim x,y\in\Omega^+}w_{xy}u^-(x)u^+(y)\nonumber\\
&    &+a\left(\frac{st-1}{2}+\frac{1-s^{p}}{2p}+\frac{1-t^{p}}{2p}\right)\sum_{x\in\Omega^+}
          \sum\limits_{y\thicksim x,y\in\Omega^-}w_{xy}u^-(y)u^+(x)\nonumber\\
&    &+b\left(\frac{1-s^2t^2}{4}-\frac{1-s^{p}}{2p}-\frac{1-t^{p}}{2p}\right)\left(\sum_{x\in\Omega^-}
          \sum\limits_{y\thicksim x,y\in\Omega^+}w_{xy}u^-(x)u^+(y)\right)^2\nonumber\\
&    &+b\left(\frac{1-s^2t^2}{4}-\frac{1-s^{p}}{2p}-\frac{1-t^{p}}{2p}\right)\left(\sum_{x\in\Omega^+}
           \sum\limits_{y\thicksim x,y\in\Omega^-}w_{xy}u^-(y)u^+(x)\right)^2\nonumber\\
&    &+b\left(\frac{1-s^2t^2}{2}-\frac{1-s^{p}}{p}-\frac{1-t^{p}}{p}\right)\sum_{x\in\Omega^-}
           \sum\limits_{y\thicksim x,y\in\Omega^+}w_{xy}u^-(x)u^+(y)
           \sum_{x\in\Omega^+}\sum\limits_{y\thicksim x,y\in\Omega^-}w_{xy}u^-(y)u^+(x)\nonumber\\
&    &+b\left[\frac{s^3t-1}{2}+\frac{3(1-s^{p})}{2p}+\frac{1-t^{p}}{2p}\right]\|u^+\|^2
           \sum_{x\in\Omega^-}\sum\limits_{y\thicksim x,y\in\Omega^+}w_{xy}u^-(x)u^+(y)\nonumber\\
&    &+b\left[\frac{s^3t-1}{2}+\frac{3(1-s^{p})}{2p}+\frac{1-t^{p}}{2p}\right]\|u^+\|^2
           \sum_{x\in\Omega^+}\sum\limits_{y\thicksim x,y\in\Omega^-}w_{xy}u^-(y)u^+(x)\nonumber\\
&    &+b\left[\frac{st^3-1}{2}+\frac{1-s^{p}}{2p}+\frac{3(1-t^{p})}{2p}\right]\|u^-\|^2
           \sum_{x\in\Omega^-}\sum\limits_{y\thicksim x,y\in\Omega^+}w_{xy}u^-(x)u^+(y)\nonumber\\
&    &+b\left[\frac{st^3-1}{2}+\frac{1-s^{p}}{2p}+\frac{3(1-t^{p})}{2p}\right]\|u^-\|^2
           \sum_{x\in\Omega^+}\sum\limits_{y\thicksim x,y\in\Omega^-}w_{xy}u^-(y)u^+(x)\nonumber\\
& =  &\frac{1-s^{p}}{p}\langle I'(u),u^+\rangle+\frac{1-t^{p}}{p}\langle I'(u),u^-\rangle
        +a\left(\frac{1-s^2}{2}-\frac{1-s^{p}}{p}\right)\|u^+\|^2
        +a\left(\frac{1-t^2}{2}-\frac{1-t^{p}}{p}\right)\|u^-\|^2\nonumber\\
&    &+b\left(\frac{1-s^4}{4}-\frac{1-s^{p}}{p}\right)\|u^+\|^4
        +b\left(\frac{1-t^4}{4}-\frac{1-t^{p}}{p}\right)\|u^-\|^4
        +\lambda\left(\frac{1-s^{\frac{2k}{m}}}{\frac{2k}{m}}-\frac{1-s^{p}}{p}\right)\int_\Omega g(x)(u^+)^{\frac{2k}{m}}d\mu\nonumber\\
&    &+\lambda\left(\frac{1-t^{\frac{2k}{m}}}{\frac{2k}{m}}-\frac{1-t^{p}}{p}\right)\int_\Omega g(x)(u^-)^{\frac{2k}{m}}d\mu
      +\frac{r(1-s^{p})+ps^{p}\ln s^r}{p^2}\int_\Omega Q(x)|u^+|^{p}d\mu\nonumber\\
&    &+\frac{r(1-t^{p})+pt^{p}\ln t^r}{p^2}\int_\Omega Q(x)|u^-|^{p}d\mu
        +b\left(\frac{1-s^2t^2}{2}-\frac{1-s^{4}}{4}-\frac{1-t^{4}}{4}\right)\|u^+\|^2\|u^-\|^2\nonumber\\
&    &+a\left(\frac{st-1}{2}+\frac{1-s^{2}}{4}+\frac{1-t^{2}}{4}\right)\sum_{x\in\Omega^-}
          \sum\limits_{y\thicksim x,y\in\Omega^+}w_{xy}u^-(x)u^+(y)\nonumber\\
&    &+a\left(\frac{st-1}{2}+\frac{1-s^{2}}{4}+\frac{1-t^{2}}{4}\right)\sum_{x\in\Omega^+}
          \sum\limits_{y\thicksim x,y\in\Omega^-}w_{xy}u^-(y)u^+(x)\nonumber\\
&    &+b\left(\frac{1-s^2t^2}{4}-\frac{1-s^{4}}{8}-\frac{1-t^{4}}{8}\right)\left(\sum_{x\in\Omega^-}
          \sum\limits_{y\thicksim x,y\in\Omega^+}w_{xy}u^-(x)u^+(y)\right)^2\nonumber\\
&    &+b\left(\frac{1-s^2t^2}{4}-\frac{1-s^{4}}{8}-\frac{1-t^{4}}{8}\right)\left(\sum_{x\in\Omega^+}
           \sum\limits_{y\thicksim x,y\in\Omega^-}w_{xy}u^-(y)u^+(x)\right)^2\nonumber\\
&    &+b\left(\frac{1-s^2t^2}{2}-\frac{1-s^{4}}{4}-\frac{1-t^{4}}{4}\right)\sum_{x\in\Omega^-}
           \sum\limits_{y\thicksim x,y\in\Omega^+}w_{xy}u^-(x)u^+(y)
           \sum_{x\in\Omega^+}\sum\limits_{y\thicksim x,y\in\Omega^-}w_{xy}u^-(y)u^+(x)\nonumber\\
&    &+b\left[\frac{s^3t-1}{2}+\frac{3(1-s^{4})}{8}+\frac{1-t^{4}}{8}\right]\|u^+\|^2
           \sum_{x\in\Omega^-}\sum\limits_{y\thicksim x,y\in\Omega^+}w_{xy}u^-(x)u^+(y)\nonumber\\
&    &+b\left[\frac{s^3t-1}{2}+\frac{3(1-s^{4})}{8}+\frac{1-t^{4}}{8}\right]\|u^+\|^2
           \sum_{x\in\Omega^+}\sum\limits_{y\thicksim x,y\in\Omega^-}w_{xy}u^-(y)u^+(x)\nonumber\\
&    &+b\left[\frac{st^3-1}{2}+\frac{1-s^{4}}{8}+\frac{3(1-t^{4})}{8}\right]\|u^-\|^2
           \sum_{x\in\Omega^-}\sum\limits_{y\thicksim x,y\in\Omega^+}w_{xy}u^-(x)u^+(y)\nonumber\\
&    &+b\left[\frac{st^3-1}{2}+\frac{1-s^{4}}{8}+\frac{3(1-t^{4})}{8}\right]\|u^-\|^2
           \sum_{x\in\Omega^+}\sum\limits_{y\thicksim x,y\in\Omega^-}w_{xy}u^-(y)u^+(x)\nonumber\\
&    &+b\left(\frac{1-s^{4}}{4}-\frac{1-s^{p}}{p}+\frac{1-t^{4}}{4}-\frac{1-t^{p}}{p}\right)\|u^+\|^2\|u^-\|^2\nonumber\\
&    &+a\left(\frac{1-s^{p}}{2p}-\frac{1-s^{2}}{4}+\frac{1-t^{p}}{2p}-\frac{1-t^{2}}{4}\right)\sum_{x\in\Omega^-}
          \sum\limits_{y\thicksim x,y\in\Omega^+}w_{xy}u^-(x)u^+(y)\nonumber\\
&    &+a\left(\frac{1-s^{p}}{2p}-\frac{1-s^{2}}{4}+\frac{1-t^{p}}{2p}-\frac{1-t^{2}}{4}\right)\sum_{x\in\Omega^+}
          \sum\limits_{y\thicksim x,y\in\Omega^-}w_{xy}u^-(y)u^+(x)\nonumber\\
&    &+b\left(\frac{1-s^{4}}{8}-\frac{1-s^{p}}{2p}+\frac{1-t^{4}}{8}-\frac{1-t^{p}}{2p}\right)\left(\sum_{x\in\Omega^-}
          \sum\limits_{y\thicksim x,y\in\Omega^+}w_{xy}u^-(x)u^+(y)\right)^2\nonumber\\
&    &+b\left(\frac{1-s^{4}}{8}-\frac{1-s^{p}}{2p}+\frac{1-t^{4}}{8}-\frac{1-t^{p}}{2p}\right)\left(\sum_{x\in\Omega^+}
           \sum\limits_{y\thicksim x,y\in\Omega^-}w_{xy}u^-(y)u^+(x)\right)^2\nonumber\\
&    &+b\left(\frac{1-s^{4}}{4}-\frac{1-s^{p}}{p}+\frac{1-t^{4}}{4}-\frac{1-t^{p}}{p}\right)\sum_{x\in\Omega^-}
           \sum\limits_{y\thicksim x,y\in\Omega^+}w_{xy}u^-(x)u^+(y)
           \sum_{x\in\Omega^+}\sum\limits_{y\thicksim x,y\in\Omega^-}w_{xy}u^-(y)u^+(x)\nonumber\\
&    &+b\left[\frac{3(1-s^{p})}{2p}-\frac{3(1-s^{4})}{8}+\frac{1-t^{p}}{2p}-\frac{1-t^{4}}{8}\right]\|u^+\|^2
           \sum_{x\in\Omega^-}\sum\limits_{y\thicksim x,y\in\Omega^+}w_{xy}u^-(x)u^+(y)\nonumber\\
&    &+b\left[\frac{3(1-s^{p})}{2p}-\frac{3(1-s^{4})}{8}+\frac{1-t^{p}}{2p}-\frac{1-t^{4}}{8}\right]\|u^+\|^2
           \sum_{x\in\Omega^+}\sum\limits_{y\thicksim x,y\in\Omega^-}w_{xy}u^-(y)u^+(x)\nonumber\\
&    &+b\left[\frac{1-s^{p}}{2p}-\frac{1-s^{4}}{8}+\frac{3(1-t^{p})}{2p}-\frac{3(1-t^{4})}{8}\right]\|u^-\|^2
           \sum_{x\in\Omega^-}\sum\limits_{y\thicksim x,y\in\Omega^+}w_{xy}u^-(x)u^+(y)\nonumber\\
&    &+b\left[\frac{1-s^{p}}{2p}-\frac{1-s^{4}}{8}+\frac{3(1-t^{p})}{2p}-\frac{3(1-t^{4})}{8}\right]\|u^-\|^2
           \sum_{x\in\Omega^+}\sum\limits_{y\thicksim x,y\in\Omega^-}w_{xy}u^-(y)u^+(x)\nonumber\\
&\geq&\frac{1-s^{p}}{p}\langle I'(u),u^+\rangle+\frac{1-t^{p}}{p}\langle I'(u),u^-\rangle
        +a\left(\frac{1-s^2}{2}-\frac{1-s^{p}}{p}\right)\|u^+\|^2
        +a\left(\frac{1-t^2}{2}-\frac{1-t^{p}}{p}\right)\|u^-\|^2\nonumber\\
&    &+b\left(\frac{1-s^4}{4}-\frac{1-s^{p}}{p}\right)\|u^+\|^4
        +b\left(\frac{1-t^4}{4}-\frac{1-t^{p}}{p}\right)\|u^-\|^4
        +\frac{b(s^2-t^2)^2}{4}\|u^+\|^2\|u^-\|^2\nonumber\\
&    &-\frac{a(s-t)^2}{4}\sum_{x\in\Omega^-}\sum\limits_{y\thicksim x,y\in\Omega^+}w_{xy}u^-(x)u^+(y)
          -\frac{a(s-t)^2}{4}\sum_{x\in\Omega^+}\sum\limits_{y\thicksim x,y\in\Omega^-}w_{xy}u^-(y)u^+(x)\nonumber\\
&    &+\frac{b(s^2-t^2)^2}{8}\left(\sum_{x\in\Omega^-}\sum\limits_{y\thicksim x,y\in\Omega^+}w_{xy}u^-(x)u^+(y)\right)^2
   +\frac{b(s^2-t^2)^2}{8}\left(\sum_{x\in\Omega^+}\sum\limits_{y\thicksim x,y\in\Omega^-}w_{xy}u^-(y)u^+(x)\right)^2\nonumber\\
&    &+\frac{b(s^2-t^2)^2}{4}\sum_{x\in\Omega^-}\sum\limits_{y\thicksim x,y\in\Omega^+}w_{xy}u^-(x)u^+(y)
           \sum_{x\in\Omega^+}\sum\limits_{y\thicksim x,y\in\Omega^-}w_{xy}u^-(y)u^+(x).
\end{eqnarray}
In the final step, we have made use of the conclusion that $\frac{s^3t-1}{2}+\frac{3(1-s^{4})}{8}+\frac{1-t^{4}}{8}<0$ which is obtained by using the Young inequality and Appendix 2. Hence, (\ref{Eq4}) holds for all $u\in \mathcal{H}^{1,2}_0(\Omega)$,\;$s,t\geq0$.
\qed
\vskip2mm
\noindent
{\bf Remark 3.1.} { By  letting $s=t$ in (\ref{Eq4}), for all $u\in \mathcal{H}^{1,2}_0(\Omega)$ and $t\geq0$, there holds}
\begin{small}
\begin{eqnarray*}
     I(u)
\geq I(tu)+\frac{1-t^{p}}{p}\langle I'(u),u\rangle
       +a\left(\frac{1-t^2}{2}-\frac{1-t^{p}}{p}\right)\left(\|u^+\|^2+\|u^-\|^2\right)
       +b\left(\frac{1-t^4}{4}-\frac{1-t^{p}}{p}\right)\left(\|u^+\|^4+\|u^-\|^4\right).
\end{eqnarray*}
\end{small}We can also obtain  the following result which is similar to that in \cite{Wen 2019} by using the similar proof as in Lemma 3.1.
\vskip2mm
\noindent
{\bf Corollary 3.2.} {\it For all $u\in \mathcal{H}^{1,2}_0(\Omega)$ and $t\geq0$, there holds}
\begin{eqnarray}
\label{Eq8}
         I(u)
& \geq & I(tu)+\frac{1-t^{p}}{p}\langle I'(u),u\rangle+a\left(\frac{1-t^2}{2}-\frac{1-t^p}{p}\right)\|u\|^2
          +b\left(\frac{1-t^4}{4}-\frac{1-t^p}{p}\right)\|u\|^4.
\end{eqnarray}
{\bf Proof.} As in Lemma 3.1, we consider the case $u\neq0$. For $u\in \mathcal{H}^{1,2}_0(\Omega)\backslash\{0\}$ and all $t\ge 0$, we have
\begin{eqnarray}
\label{Equ1}
\int_\Omega |tu|^{p}\ln |tu|^r d\mu=\int_\Omega \left(|tu|^{p}\ln |u|^r+|tu|^p\ln t^r\right) d\mu.
\end{eqnarray}
Because the function $f(x)=\frac{1-a^x}{x}$ is monotonically decreasing on $(0,+\infty)$ for $a>0$ and $a\neq 1$ \big(see Appendix 2\big), in combination with (\ref{eq9}),\;(\ref{eq12}),
(\ref{Eq5}) and (\ref{Equ1}), we have
\begin{eqnarray}
\label{Equ2}
&    &I(u)-I(tu)\nonumber\\
& =  &\frac{a}{2}\left(\|u\|^2-\|tu\|^2\right)+\frac{b}{4}\left(\|u\|^4-\|tu\|^4\right)
      +\frac{m}{2k}\int_\Omega\lambda g(x)\left[u^{\frac{2k}{m}}-(tu)^{\frac{2k}{m}}\right]d\mu
      +\frac{r}{p^2}\int_\Omega Q(x)\left(|u|^p-|tu|^p\right)d\mu\nonumber\\
&    &-\frac{1}{p}\int_\Omega Q(x)\Big(|u|^{p}\ln |u|^r-|tu|^{p}\ln |tu|^r\Big)d\mu\nonumber\\
& =  &\frac{a(1-t^2)}{2}\|u\|^2+\frac{b(1-t^4)}{4}\|u\|^4
         +\frac{1-t^{\frac{2k}{m}}}{\frac{2k}{m}}\int_\Omega\lambda g(x)u^{\frac{2k}{m}}d\mu
         +\frac{r(1-t^p)}{p^2}\int_\Omega Q(x)|u|^pd\mu\nonumber\\
&    &-\frac{1-t^p}{p}\int_\Omega Q(x)|u|^p\ln{|u|^r}d\mu+\frac{1}{p}\int_\Omega Q(x)|tu|^p\ln {t^r}d\mu\nonumber\\
& =  &\frac{1-t^p}{p}\left[\langle I'(u),u\rangle-a\|u\|^2-b\|u\|^4-\int_\Omega\lambda g(x)u^{\frac{2k}{m}}d\mu\right]
      +\frac{a(1-t^2)}{2}\|u\|^2+\frac{b(1-t^4)}{4}\|u\|^4\nonumber\\
&    &+\frac{1-t^{\frac{2k}{m}}}{\frac{2k}{m}}\int_\Omega\lambda g(x)u^{\frac{2k}{m}}d\mu
      +\frac{r(1-t^p)}{p^2}\int_\Omega Q(x)|u|^pd\mu
      +\frac{1}{p}\int_\Omega Q(x)|tu|^p\ln {t^r}d\mu\nonumber\\
& =  &\frac{1-t^p}{p}\langle I'(u),u\rangle+a\left(\frac{1-t^2}{2}-\frac{1-t^p}{p}\right)\|u\|^2
       +b\left(\frac{1-t^4}{4}-\frac{1-t^p}{p}\right)\|u\|^4\nonumber\\
&    &+\left(\frac{1-t^{\frac{2k}{m}}}{\frac{2k}{m}}-\frac{1-t^p}{p}\right)\int_\Omega\lambda g(x)u^{\frac{2k}{m}}d\mu
       +\frac{r(1-t^p)+pt^p\ln t^r}{p^2}\int_\Omega Q(x)|u|^pd\mu\nonumber\\
&\geq&\frac{1-t^p}{p}\langle I'(u),u\rangle+a\left(\frac{1-t^2}{2}-\frac{1-t^p}{p}\right)\|u\|^2
        +b\left(\frac{1-t^4}{4}-\frac{1-t^p}{p}\right)\|u\|^4.
\end{eqnarray}
Hence, (\ref{Eq8}) holds for all $u\in \mathcal{H}^{1,2}_0(\Omega)$,\;$t\geq0$. \qed
\par
Note that $p>4$, $u^-(y)u^+(x)<0$, $u^+(y)u^-(x)<0$ and the function $f(x)=\frac{1-a^x}{x}$ is monotonically decreasing on $(0,+\infty)$ for $a>0$ and $a\neq 1$. Then in view of Lemma 3.1, we have the following corollary.
\vskip2mm
\noindent
{\bf Corollary 3.3.} {\it For any $u\in \mathcal{M}$, there exists $I(u)=\max_{s,t\geq0}I(su^++tu^-)$.}
\par
In view of Corollary 3.2 or Remark 3.1, we have the following corollary.\\
\noindent
{\bf Corollary 3.4.} {\it For any $u\in \mathcal{N}$, there exists $I(u)=\max_{t\geq0}I(tu)$.}
\vskip2mm
\noindent
{\bf Lemma 3.5.} {\it For any $u\in \mathcal{H}^{1,2}_0(\Omega)\setminus\{0\}$, there exists a unique $t_0>0$ such that $t_0u\in\mathcal{N}$.}\\
{\bf Proof.} First, we prove that the existence of $t_0$. For any $u\in \mathcal{H}^{1,2}_0(\Omega)\setminus\{0\}$, let $u\in\mathcal{N}$ be fixed and define a function $f(t)=\langle I'(tu),tu\rangle$ on $(0,+\infty)$. Note that $\lambda\ge 0$ and $g(x)>0$ for all $x\in \Omega$. Then, according to (\ref{eq11}) and Lemma 2.1, there exist positive constants $\varepsilon_1<\frac{a}{K_2^2}$ and $C_{\varepsilon_1}$ such that
\begin{eqnarray}
\label{a1}
         f(t)
&  =   & at^2\|u\|^2+bt^4\|u\|^4-\int_\Omega Q(x)|tu|^{p}\ln|tu|^rd\mu
        +\lambda t^{\frac{2k}{m}}\int_\Omega g(x)u^{\frac{2k}{m}}d\mu\nonumber\\
& \geq & at^2\|u\|^2+bt^4\|u\|^4-\int_\Omega \varepsilon_1|tu|^2d\mu
        -\int_\Omega C_{\varepsilon_1}|tu|^ld\mu+\lambda t^{\frac{2k}{m}}\int_\Omega g(x)u^{\frac{2k}{m}}d\mu\nonumber\\
& \geq & at^2\|u\|^2-K_2^2\varepsilon_1 t^2\|u\|^2+bt^4\|u\|^4
        -t^l\int_\Omega C_{\varepsilon_1}|u|^ld\mu+\lambda t^{\frac{2k}{m}}\int_\Omega g(x)u^{\frac{2k}{m}}d\mu.
\end{eqnarray}
Then it follows from $l>p>4$ that $f(t)>0$ for all sufficiently small $t$.
\par
One the other hand, noting that  $Q(x)>0$ for all $x\in \Omega$, we have
\begin{eqnarray}
\label{a2}
         f(t)
&  =  & at^2\|u\|^2+bt^4\|u\|^4-t^p\ln t^r\int_\Omega Q(x)|u|^{p}d\mu
        -t^p\int_\Omega Q(x)|u|^{p}\ln|u|^rd\mu+\lambda t^{\frac{2k}{m}}\int_\Omega g(x)u^{\frac{2k}{m}}d\mu\nonumber\\
&  \le & at^2\|u\|^2+bt^4\|u\|^4-t^p\ln t^r\int_\Omega Q(x)|u|^{p}d\mu
        +t^p\int_\Omega Q(x)|u|^{p}|\ln|u|^r|d\mu+\lambda t^{\frac{2k}{m}}\int_\Omega g(x)u^{\frac{2k}{m}}d\mu.\nonumber\\
\end{eqnarray}
Then it follows from $\frac{2k}{m}\leq p$ and $p>4$  that there exists large enough $t$ such that $f(t)<0$.
Thus by the continuity of $f(t)$, there exists $t_0>0$ such that $f(t_0)=0$, which implies that there appears a number $t_0>0$ such that $t_0u\in\mathcal{N}$.
\par
Next, we prove the uniqueness of $t_0$. Arguing by contradiction, we suppose that there are $u\in \mathcal{H}^{1,2}_0(\Omega)\backslash\{0\}$ and two positive constants $t_1\neq t_2$ such that $t_1u\in\mathcal{N}$ and $t_2u\in\mathcal{N}$. Note that the function $f(x)=\frac{1-a^x}{x}$ is strictly monotonically decreasing on $(0,+\infty)$ for $a>0$ and $a\neq 1$. In Corollary 3.2, taking $u:=t_1u$ and  $t:=\frac{t_2}{t_1}$, there holds
\begin{eqnarray}\label{a11}
      I(t_1u)
 \geq I(t_2u)+at_1^2\left[\frac{1-(\frac{t_2}{t_1})^2}{2}-\frac{1-(\frac{t_2}{t_1})^{p}}{p}\right]\|u\|^2
        +bt_1^4\left[\frac{1-(\frac{t_2}{t_1})^4}{4}-\frac{1-(\frac{t_2}{t_1})^{p}}{p}\right]\|u\|^4
   >  I(t_2u),
\end{eqnarray}
and taking $u:=t_2u$ and  $t:=\frac{t_1}{t_2}$, there holds
\begin{eqnarray}\label{a22}
      I(t_2u)
 \geq I(t_1u)+at_2^2\left[\frac{1-(\frac{t_1}{t_2})^2}{2}-\frac{1-(\frac{t_1}{t_2})^{p}}{p}\right]\|u\|^2
        +bt_2^4\left[\frac{1-(\frac{t_1}{t_2})^4}{4}-\frac{1-(\frac{t_1}{t_2})^{p}}{p}\right]\|u\|^4
   >  I(t_1u),
\end{eqnarray}
which is a contradiction. So $t_1=t_2$, that is, there is a unique $t_0>0$ such that $t_0u\in\mathcal{N}$.
\qed

\vskip2mm
\noindent
{\bf Lemma 3.6.} {\it For any $u\in \mathcal{H}^{1,2}_0(\Omega)$ with $u^\pm\neq 0$, there exists a unique pair of positive numbers $(s_0,t_0)$ such that $s_0u^++t_0u^-\in\mathcal{M}$.}\\
{\bf Proof.} First, for any $u\in \mathcal{H}^{1,2}_0(\Omega)$ with $u^\pm\neq 0$, we prove the existence of $(s_0,t_0)$. By Proposition 3.2, we let
\begin{eqnarray}
\label{G1}
     G(s,t)
& = &\langle I'(su^++tu^-),su^+\rangle\nonumber\\
& = &as^2\|u^+\|^2+bs^4\|u^+\|^4-\int_\Omega Q(x)|su^+|^p\ln{|su^+|^r}d\mu
      +\int_\Omega\lambda g(x)(su^+)^{\frac{2k}{m}}d\mu\nonumber\\
&   &-\frac{ast}{2}\sum_{x\in\Omega^-}\sum\limits_{y\thicksim x,y\in\Omega^+}w_{xy}u^-(x)u^+(y)
      -\frac{ast}{2}\sum_{x\in\Omega^+}\sum\limits_{y\thicksim x,y\in\Omega^-}w_{xy}u^-(y)u^+(x)\nonumber\\
&   &+\frac{bs^2t^2}{2}\left(\sum_{x\in\Omega^-}\sum\limits_{y\thicksim x,y\in\Omega^+}w_{xy}u^-(x)u^+(y)\right)^2
      +\frac{bs^2t^2}{2}\left(\sum_{x\in\Omega^+}\sum\limits_{y\thicksim x,y\in\Omega^-}w_{xy}u^-(y)u^+(x)\right)^2\nonumber\\
&   &+bs^2t^2\|u^+\|^2\|u^-\|^2+bs^2t^2\sum_{x\in\Omega^+}\sum\limits_{y\thicksim x,y\in\Omega^-}w_{xy}u^-(y)u^+(x)
           \sum_{x\in\Omega^-}\sum\limits_{y\thicksim x,y\in\Omega^+}w_{xy}u^-(x)u^+(y)\nonumber\\
&   &-\frac{3bs^3t}{2}\|u^+\|^2\sum_{x\in\Omega^-}\sum\limits_{y\thicksim x,y\in\Omega^+}w_{xy}u^-(x)u^+(y)
        -\frac{3bs^3t}{2}\|u^+\|^2\sum_{x\in\Omega^+}\sum\limits_{y\thicksim x,y\in\Omega^-}w_{xy}u^-(y)u^+(x)\nonumber\\
&   &-\frac{bst^3}{2}\|u^-\|^2\sum_{x\in\Omega^-}\sum\limits_{y\thicksim x,y\in\Omega^+}w_{xy}u^-(x)u^+(y)
        -\frac{bst^3}{2}\|u^-\|^2\sum_{x\in\Omega^+}\sum\limits_{y\thicksim x,y\in\Omega^-}w_{xy}u^-(y)u^+(x)
\end{eqnarray}
and
\begin{eqnarray}
\label{G2}
     H(s,t)
& = &\langle I'(su^++tu^-),tu^-\rangle\nonumber\\
& = &at^2\|u^-\|^2+bt^4\|u^-\|^4-\int_\Omega Q(x)|tu^-|^p\ln{|tu^-|^r}d\mu
     +\int_\Omega\lambda g(x)(tu^-)^{\frac{2k}{m}}d\mu\nonumber\\
&   &-\frac{ast}{2}\sum_{x\in\Omega^-}\sum\limits_{y\thicksim x,y\in\Omega^+}w_{xy}u^-(x)u^+(y)
     -\frac{ast}{2}\sum_{x\in\Omega^+}\sum\limits_{y\thicksim x,y\in\Omega^-}w_{xy}u^-(y)u^+(x)\nonumber\\
&   &+\frac{bs^2t^2}{2}\left(\sum_{x\in\Omega^-}\sum\limits_{y\thicksim x,y\in\Omega^+}w_{xy}u^-(x)u^+(y)\right)^2
      +\frac{bs^2t^2}{2}\left(\sum_{x\in\Omega^+}\sum\limits_{y\thicksim x,y\in\Omega^-}w_{xy}u^-(y)u^+(x)\right)^2\nonumber\\
&   &+bs^2t^2\|u^+\|^2\|u^-\|^2+bs^2t^2\sum_{x\in\Omega^+}\sum\limits_{y\thicksim x,y\in\Omega^-}w_{xy}u^-(y)u^+(x)
           \sum_{x\in\Omega^-}\sum\limits_{y\thicksim x,y\in\Omega^+}w_{xy}u^-(x)u^+(y)\nonumber\\
&   &-\frac{3bst^3}{2}\|u^-\|^2\sum_{x\in\Omega^-}\sum\limits_{y\thicksim x,y\in\Omega^+}w_{xy}u^-(x)u^+(y)
        -\frac{3bst^3}{2}\|u^-\|^2\sum_{x\in\Omega^+}\sum\limits_{y\thicksim x,y\in\Omega^-}w_{xy}u^-(y)u^+(x)\nonumber\\
&   &-\frac{bs^3t}{2}\|u^+\|^2\sum_{x\in\Omega^-}\sum\limits_{y\thicksim x,y\in\Omega^+}w_{xy}u^-(x)u^+(y)
        -\frac{bs^3t}{2}\|u^+\|^2\sum_{x\in\Omega^+}\sum\limits_{y\thicksim x,y\in\Omega^-}w_{xy}u^-(y)u^+(x).
\end{eqnarray}
Note that $p>4$, $u^-(y)u^+(x)<0$ and $u^+(y)u^-(x)<0$. Similar to the argument of (\ref{a1}) and (\ref{a2}), it is easy to see that $G(s,s)>0$ and $H(s,s)>0$ for $s>0$ sufficiently small and $G(t,t)<0$ and $H(t,t)<0$ for $t>0$ large enough. Thus, there exist constants $0<\alpha<\beta$ such that
\begin{eqnarray}
\label{G3}
G(\alpha,\alpha)>0,\;H(\alpha,\alpha)>0, \;G(\beta,\beta)<0,\;H(\beta,\beta)<0.
\end{eqnarray}
From (\ref{G1}), (\ref{G2}) and (\ref{G3}), $\forall s,t\in [\alpha,\beta]$, we have
\begin{eqnarray}
\label{G4}
         G(\alpha,t)
&\geq&a\alpha^2\|u^+\|^2+b\alpha^4\|u^+\|^4-\int_\Omega Q(x)|\alpha u^+|^p\ln{|\alpha u^+|^r}d\mu
         +\int_\Omega\lambda g(x)(\alpha u^+)^{\frac{2k}{m}}d\mu\nonumber\\
&    &-\frac{\alpha^2}{2}\sum_{x\in\Omega^-}\sum\limits_{y\thicksim x,y\in\Omega^+}w_{xy}u^-(x)u^+(y)
         -\frac{\alpha^2}{2}\sum_{x\in\Omega^+}\sum\limits_{y\thicksim x,y\in\Omega^-}w_{xy}u^-(y)u^+(x)\nonumber\\
&    &+\frac{b\alpha^4}{2}\left(\sum_{x\in\Omega^-}\sum\limits_{y\thicksim x,y\in\Omega^+}w_{xy}u^-(x)u^+(y)\right)^2
        +\frac{b\alpha^4}{2}\left(\sum_{x\in\Omega^+}\sum\limits_{y\thicksim x,y\in\Omega^-}w_{xy}u^-(y)u^+(x)\right)^2\nonumber\\
&    &+b\alpha^4\|u^+\|^2\|u^-\|^2+b\alpha^4\sum_{x\in\Omega^+}\sum\limits_{y\thicksim x,y\in\Omega^-}w_{xy}u^-(y)u^+(x)
           \sum_{x\in\Omega^-}\sum\limits_{y\thicksim x,y\in\Omega^+}w_{xy}u^-(x)u^+(y)\nonumber\\
&    &-\frac{3b\alpha^4}{2}\|u^+\|^2\sum_{x\in\Omega^-}\sum\limits_{y\thicksim x,y\in\Omega^+}w_{xy}u^-(x)u^+(y)
        -\frac{3b\alpha^4}{2}\|u^+\|^2\sum_{x\in\Omega^+}\sum\limits_{y\thicksim x,y\in\Omega^-}w_{xy}u^-(y)u^+(x)\nonumber\\
&    &-\frac{b\alpha^4}{2}\|u^-\|^2\sum_{x\in\Omega^-}\sum\limits_{y\thicksim x,y\in\Omega^+}w_{xy}u^-(x)u^+(y)
        -\frac{b\alpha^4}{2}\|u^-\|^2\sum_{x\in\Omega^+}\sum\limits_{y\thicksim x,y\in\Omega^-}w_{xy}u^-(y)u^+(x)\nonumber\\
&  = &G(\alpha,\alpha)\nonumber\\
&  > &0,
\end{eqnarray}
and using the same method as (\ref{G4}), we can also obtain that
\begin{eqnarray}
\label{G5}
G(\beta,t)\leq G(\beta,\beta)<0,\;\;H(s,\alpha)\geq H(\alpha,\alpha)>0\;\;\text{and}\;\;H(s,\beta)\leq H(\beta,\beta)<0.
\end{eqnarray}
So, combining (\ref{G4}), (\ref{G5}), with the Poincar\'{e}-Miranda Theorem \cite{Miranda 1940}, it is easy to see that there exists a point $(s_0,t_0)$ with $\alpha<s_0,t_0<\beta$ such that $G(s_0,t_0)=H(s_0,t_0)=0$, which implies that there exists a pair of positive numbers $(s_0,t_0)$ such that $s_0u^++t_0u^-\in\mathcal{M}$.
\par
Next, we prove the uniqueness of $(s_0,t_0)$. Arguing by contradiction, we assume that there exist two pairs of positive numbers $(s_1,t_1)$ and $(s_2,t_2)$ such that $(s_1,t_1)\not=(s_2,t_2)$, $s_1u^++t_1u^-\in\mathcal{M}$ and $s_2u^++t_2u^-\in\mathcal{M}$. Note that $p>4$ and $f(x)=\frac{1-a^x}{x}$ is strictly monotonically decreasing on $(0,+\infty)$ for $a>0$ and $a\neq1$. In Lemma 3.1, taking $u:=s_1u^++t_1u^-$, $s:=\frac{s_2}{s_1}$ and  $t:=\frac{t_2}{t_1}$, respectively, there holds
\begin{eqnarray}
\label{s1}
         I(s_1u^++t_1u^-)
&\geq&I(s_2u^++t_2u^-)+as_1^2\left[\frac{1-(\frac{s_2}{s_1})^2}{2}-\frac{1-(\frac{s_2}{s_1})^{p}}{p}\right]\|u^+\|^2\nonumber\\
&    &+at_1^2\left[\frac{1-(\frac{t_2}{t_1})^2}{2}-\frac{1-(\frac{t_2}{t_1})^{p}}{p}\right]\|u^-\|^2
         +bs_1^4\left[\frac{1-(\frac{s_2}{s_1})^4}{4}-\frac{1-(\frac{s_2}{s_1})^{p}}{p}\right]\|u^+\|^4\nonumber\\
&    &+bt_1^4\left[\frac{1-(\frac{t_2}{t_1})^4}{4}-\frac{1-(\frac{t_2}{t_1})^{p}}{p}\right]\|u^-\|^4
        -\frac{as_1t_1(\frac{s_2}{s_1}-\frac{t_2}{t_1})^2}{4}\sum_{x\in\Omega^-}\sum\limits_{y\thicksim x,y\in\Omega^+}w_{xy}u^-(x)u^+(y)\nonumber\\
&    &-\frac{as_1t_1(\frac{s_2}{s_1}-\frac{t_2}{t_1})^2}{4}\sum_{x\in\Omega^+}\sum\limits_{y\thicksim x,y\in\Omega^-}w_{xy}u^-(y)u^+(x)
       +\frac{bs_1^2t_1^2[(\frac{s_2}{s_1})^2-(\frac{t_2}{t_1})^2]^2}{4}\|u^+\|^2\|u^-\|^2\nonumber\\
&    &+\frac{bs_1^2t_1^2[(\frac{s_2}{s_1})^2-(\frac{t_2}{t_1})^2]^2}{4}\sum_{x\in\Omega^-}\sum\limits_{y\thicksim x,y\in\Omega^+}w_{xy}u^-(x)u^+(y)
           \sum_{x\in\Omega^+}\sum\limits_{y\thicksim x,y\in\Omega^-}w_{xy}u^-(y)u^+(x)\nonumber\\
&    &+\frac{bs_1^2t_1^2[(\frac{s_2}{s_1})^2-(\frac{t_2}{t_1})^2]^2}{8}\left(\sum_{x\in\Omega^-}\sum\limits_{y\thicksim x,y\in\Omega^+}w_{xy}u^-(x)u^+(y)\right)^2\nonumber\\
&    &+\frac{bs_1^2t_1^2[(\frac{s_2}{s_1})^2-(\frac{t_2}{t_1})^2]^2}{8}\left(\sum_{x\in\Omega^+}\sum\limits_{y\thicksim x,y\in\Omega^-}w_{xy}u^-(y)u^+(x)\right)^2\nonumber\\
&  > & I(s_2u^++t_2u^-),
\end{eqnarray}
and taking $u:=s_2u^++t_2u^-$, $s:=\frac{s_1}{s_2}$ and  $t:=\frac{t_1}{t_2}$, respectively, we can obtain that
\begin{eqnarray}
\label{s2}
I(s_2u^++t_2u^-)>I(s_1u^++t_1u^-).
\end{eqnarray}
Therefore, it is easy to see that (\ref{s1}) contradicts with (\ref{s2}).  Thus, $(s_1,t_1)=(s_2,t_2)$, that is, there appears a unique pair of positive numbers $(s_0,t_0)$ such that $s_0u^++t_0u^-\in\mathcal{M}$.
\qed

\vskip2mm
\noindent
{\bf Lemma 3.7.} {\it There exist the following minimax characterizations}
\begin{eqnarray}
\label{a3}
  \inf_{u\in\mathcal{N}}I(u)
= :c
= \inf_{u\in \mathcal{H}^{1,2}_0(\Omega),u\neq 0}\max_{t\geq 0}I(tu)
\end{eqnarray}
{\it and}
\begin{eqnarray}
\label{a4}
  \inf_{u\in\mathcal{M}}I(u)
= :m
= \inf_{u\in \mathcal{H}^{1,2}_0(\Omega),u^\pm\neq 0}\max_{s,t\geq 0}I(su^++tu^-).
\end{eqnarray}
\vskip0mm
\noindent
{\bf Proof.} For one thing, it follows from Corollary 3.4 and the definition of $\mathcal{N}$ that
$$
     \inf_{u\in\mathcal{N}}I(u)
  =  \inf_{u\in \mathcal{N}}\max_{t\geq 0}I(tu)
\geq \inf_{u\in \mathcal{H}^{1,2}_0(\Omega),u\neq 0}\max_{t\geq 0}I(tu).
$$
For another, by Lemma 3.5, for any $u\in \mathcal{H}^{1,2}_0(\Omega)$ and $u\neq 0$, there exists $t_u>0$ such that $t_u u\in \mathcal{N}$. Then there exists
\begin{eqnarray}
\label{a5}
     \max_{t\ge 0} I(tu)
\geq I(t_u u)
\geq \inf_{u\in\mathcal{N}}I(u).
\end{eqnarray}
Furthermore, we have
$$
     \inf_{u\in \mathcal{H}^{1,2}_0(\Omega),u\neq0}\max_{t\ge 0} I(tu)
\geq \inf_{u\in \mathcal{H}^{1,2}_0(\Omega),u\neq0}I(t_u u)
\geq \inf_{u\in\mathcal{N}}I(u).
$$
Hence, (\ref{a3}) holds. In the same way, by Corollary 3.3, the definition of $\mathcal{M}$ and Lemma 3.6,  we can obtain that (\ref{a4}) also holds.

\vskip2mm
\noindent
{\bf Lemma 3.8.} {\it $c>0$ and $m>0$ are achieved.}
\vskip0mm
\noindent
{\bf Proof.} For any $u\in\mathcal{M}$, there exists $\langle I'(u),u\rangle=0$.
According to (\ref{eq11}), (\ref{eq12}) and  Lemma 2.1, for $\varepsilon_2=\frac{a}{2K_2^2}$, there exists a positive constant $C_{\varepsilon_2}$ such that
\begin{eqnarray}
\label{Eq19}
          \|u\|^2
&  \leq  &\|u\|^2+\frac{b}{a}\|u\|^4+\frac{1}{a}\int_\Omega\lambda g(x)u^{\frac{2k}{m}}d\mu
         =\frac{1}{a}\int_{\Omega}Q(x)|u|^{p}\ln |u|^rd\mu\nonumber\\
&  \leq  &\varepsilon_2\|u\|_{L^{2}(\Omega)}^{2}+C_{\varepsilon_2}\|u\|_{L^{l}(\Omega)}^{l}\nonumber\\
&  \leq  &\frac{1}{2}\|u\|^2+C_{\varepsilon_2} K_l^l\|u\|^l.
\end{eqnarray}
Since $4<p<l$, then $\|u\|\geq \rho:=\frac{1}{\left(2C_{\varepsilon_2}K_l^l\right)^{l-2}}$ for any $u\in\mathcal{M}$.
\par
Let $\{u_n\}\subset \mathcal{M}$ be such that $I(u_n)\rightarrow m.$ By (\ref{eq9}) and (\ref{eq12}), there holds
\begin{eqnarray*}
         m+o(1)
&  =   & I(u_n)-\frac{1}{p}\langle I'(u_n),u_n\rangle\nonumber\\
&  =   & a\left(\frac{1}{2}-\frac{1}{p}\right)\|u_n\|^2+b\left(\frac{1}{4}-\frac{1}{p}\right)\|u_n\|^4
           +\frac{r}{p^2}\int_\Omega Q(x)|u_n|^pd\mu
           +\left(\frac{m}{2k}-\frac{1}{p}\right)\int_\Omega\lambda g(x)u_n^{\frac{2k}{m}}d\mu\nonumber\\
& \geq & a\left(\frac{1}{2}-\frac{1}{p}\right)\|u_n\|^2.
\end{eqnarray*}
This shows that the sequence $\{u_n\}$ is bounded in $\mathcal{H}^{1,2}_0(\Omega)$, that is, there exists a $M>0$ such that $\|u_n\|\le M$. Because $\mathcal{H}^{1,2}_0(\Omega)$ has finite dimension, then by Lemma 2.1, there exists a subsequence, still denoted by
$\{u_n\}$, and a function $\tilde{u}$, such that $u_n\rightarrow \tilde{u}$ in $\mathcal{H}^{1,2}_0(\Omega)$. Furthermore, in virtue of Lemma 2.1, we have
\begin{eqnarray}
\label{Eq20}
 \begin{cases}
 u_n^\pm \to \tilde{u}^\pm &\text{in}\;\mathcal{H}^{1,2}_0(\Omega),\\
 u_n^\pm(x)\to \tilde{u}^\pm(x) & \mbox{for each } x\in \Omega,\\
 u_n(x)\to \tilde{u}(x) & \mbox{for each } x\in \Omega,\\
 u_n^\pm\rightarrow\tilde{u}^\pm & \text{in} \; L^s(\Omega), s\in[2,+\infty).
 \end{cases}
\end{eqnarray}
Since $\{u_n\}\subset \mathcal{M}$, there exists $\langle I'(u_n),u_n^\pm\rangle=0$ and then by Proposition 3.2, we have
\begin{eqnarray}\label{aa1}
&   & a\|u_n^{\pm}\|^2+b\|u_n^{\pm}\|^4+\int_\Omega\lambda  g(x)(u_n^{\pm})^{\frac{2k}{m}}d\mu
        -\int_\Omega Q(x)|u_n^{\pm}|^{p}\ln |u_n^{\pm}|^rd\mu\nonumber\\
& = & \frac{a}{2}\sum_{x\in\Omega^-}\sum\limits_{y\thicksim x,y\in\Omega^+}w_{xy}u_n^-(x)u_n^+(y)
        +\frac{a}{2}\sum_{x\in\Omega^+}\sum\limits_{y\thicksim x,y\in\Omega^-}w_{xy}u_n^-(y)u_n^+(x)\nonumber\\
&   & -\frac{b}{2}\left(\sum_{x\in\Omega^-}\sum\limits_{y\thicksim x,y\in\Omega^+}w_{xy}u_n^-(x)u_n^+(y)\right)^2
        -\frac{b}{2}\left(\sum_{x\in\Omega^+}\sum\limits_{y\thicksim x,y\in\Omega^-}w_{xy}u_n^-(y)u_n^+(x)\right)^2\nonumber\\
&   & -b\|u_n^+\|^2\|u_n^-\|^2-b\sum_{x\in\Omega^+}\sum\limits_{y\thicksim x,y\in\Omega^-}w_{xy}u_n^-(y)u_n^+(x)
           \sum_{x\in\Omega^-}\sum\limits_{y\thicksim x,y\in\Omega^+}w_{xy}u_n^-(x)u_n^+(y)\nonumber\\
&   & +\frac{3b}{2}\|u_n^\pm\|^2\sum_{x\in\Omega^-}\sum\limits_{y\thicksim x,y\in\Omega^+}w_{xy}u_n^-(x)u_n^+(y)
         +\frac{3b}{2}\|u_n^\pm\|^2\sum_{x\in\Omega^+}\sum\limits_{y\thicksim x,y\in\Omega^-}w_{xy}u_n^-(y)u_n^+(x)\nonumber\\
&   & +\frac{b}{2}\|u_n^\mp\|^2\sum_{x\in\Omega^-}\sum\limits_{y\thicksim x,y\in\Omega^+}w_{xy}u_n^-(x)u_n^+(y)
         +\frac{b}{2}\|u_n^\mp\|^2\sum_{x\in\Omega^+}\sum\limits_{y\thicksim x,y\in\Omega^-}w_{xy}u_n^-(y)u_n^+(x).
\end{eqnarray}
By (\ref{eq11}), (\ref{Eq19}) and (\ref{aa1}), there exists $\varepsilon_3\in (0,\frac{a\rho^2}{K_2^2M^2})$ and a positive constant $C_{\varepsilon_3}$ such that
\begin{eqnarray*}
         a\rho^2
&\leq&  a\|u_n^\pm\|^2 \nonumber\\
&\leq&  a\|u_n^\pm\|^2+b\|u_n^\pm\|^4+\int_\Omega\lambda g(x)(u_n^\pm)^{\frac{2k}{m}} d\mu\nonumber\\
& =  &  \int_{\Omega}Q(x)|u_n^\pm|^{p}\ln|u_n^\pm|^rd\mu
          +\frac{a}{2}\sum_{x\in\Omega^-}\sum\limits_{y\thicksim x,y\in\Omega^+}w_{xy}u^-(x)u^+(y)
          +\frac{a}{2}\sum_{x\in\Omega^+}\sum\limits_{y\thicksim x,y\in\Omega^-}w_{xy}u^-(y)u^+(x)\nonumber\\
&    & -\frac{b}{2}\left(\sum_{x\in\Omega^-}\sum\limits_{y\thicksim x,y\in\Omega^+}w_{xy}u_n^-(x)u_n^+(y)\right)^2
          -\frac{b}{2}\left(\sum_{x\in\Omega^+}\sum\limits_{y\thicksim x,y\in\Omega^-}w_{xy}u_n^-(y)u_n^+(x)\right)^2\nonumber\\
&    & -b\|u_n^+\|^2\|u_n^-\|^2-b\sum_{x\in\Omega^+}\sum\limits_{y\thicksim x,y\in\Omega^-}w_{xy}u_n^-(y)u_n^+(x)
           \sum_{x\in\Omega^-}\sum\limits_{y\thicksim x,y\in\Omega^+}w_{xy}u_n^-(x)u_n^+(y)\nonumber\\
&    & +\frac{3b}{2}\|u_n^\pm\|^2\sum_{x\in\Omega^-}\sum\limits_{y\thicksim x,y\in\Omega^+}w_{xy}u_n^-(x)u_n^+(y)
         +\frac{3b}{2}\|u_n^\pm\|^2\sum_{x\in\Omega^+}\sum\limits_{y\thicksim x,y\in\Omega^-}w_{xy}u_n^-(y)u_n^+(x)\nonumber\\
&    & +\frac{b}{2}\|u_n^\mp\|^2\sum_{x\in\Omega^-}\sum\limits_{y\thicksim x,y\in\Omega^+}w_{xy}u_n^-(x)u_n^+(y)
         +\frac{b}{2}\|u_n^\mp\|^2\sum_{x\in\Omega^+}\sum\limits_{y\thicksim x,y\in\Omega^-}w_{xy}u_n^-(y)u_n^+(x)\nonumber\\
&\leq&  \varepsilon_3\|u_n^\pm\|_{L^{2}(\Omega)}^{2}+C_{\varepsilon_3}\|u_n^\pm\|_{L^{l}(\Omega)}^{l}\nonumber\\
&\leq&  \varepsilon_3 K_2^2\|u_n^\pm\|^{2}+C_{\varepsilon_3}\|u_n^\pm\|_{L^{l}(\Omega)}^{l}\nonumber\\
&\leq&  \varepsilon_3 K_2^2M^{2}+C_{\varepsilon_3}\|u_n^\pm\|_{L^{l}(\Omega)}^{l}.
\end{eqnarray*}
Thus, $\|u_n^\pm\|_{L^{l}(\Omega)}^{l}\geq \frac{a\rho^2-\varepsilon_3 K_2^2M^2}{C_{\varepsilon_3}}>0$, which implies that $\tilde{u}^\pm\neq 0$. Note that the definition of integral of a function $u$ is like (\ref{eq4}) and $\Omega$ is of finitely
many points.  Then according to (\ref{eq12}), (\ref{Eq20}) and (\ref{aa1}), there holds
\begin{eqnarray*}
&   & a\|\tilde{u}^\pm\|^2+b\|\tilde{u}^\pm\|^4+\int_\Omega\lambda g(x)(\tilde{u}^\pm)^{\frac{2k}{m}} d\mu\\
& = & \lim_{n\rightarrow\infty}a\|u_n^\pm\|^2+\lim_{n\rightarrow\infty}b\|u_n^\pm\|^4
        +\lim_{n\rightarrow\infty}\int_\Omega\lambda g(x)(u_n^\pm)^{\frac{2k}{m}} d\mu\\
& = & \lim_{n\rightarrow\infty}\left[\int_{\Omega}Q(x)|u_n^\pm|^{p}\ln |u_n^\pm|^rd\mu
        +\frac{a}{2}\sum_{x\in\Omega^-}\sum\limits_{y\thicksim x,y\in\Omega^+}w_{xy}u_n^-(x)u_n^+(y)
        +\frac{a}{2}\sum_{x\in\Omega^+}\sum\limits_{y\thicksim x,y\in\Omega^-}w_{xy}u_n^-(y)u_n^+(x)\right.\\
&   & \left.-\frac{b}{2}\left(\sum_{x\in\Omega^-}\sum\limits_{y\thicksim x,y\in\Omega^+}w_{xy}u_n^-(x)u_n^+(y)\right)^2
          -\frac{b}{2}\left(\sum_{x\in\Omega^+}\sum\limits_{y\thicksim x,y\in\Omega^-}w_{xy}u_n^-(y)u_n^+(x)\right)^2\right.\\
&   & \left.-b\|u_n^+\|^2\|u_n^-\|^2-b\sum_{x\in\Omega^+}\sum\limits_{y\thicksim x,y\in\Omega^-}w_{xy}u_n^-(y)u_n^+(x)
           \sum_{x\in\Omega^-}\sum\limits_{y\thicksim x,y\in\Omega^+}w_{xy}u_n^-(x)u_n^+(y)\right.\\
&   & \left.+\frac{3b}{2}\|u_n^\pm\|^2\sum_{x\in\Omega^-}\sum\limits_{y\thicksim x,y\in\Omega^+}w_{xy}u_n^-(x)u_n^+(y)
         +\frac{3b}{2}\|u_n^\pm\|^2\sum_{x\in\Omega^+}\sum\limits_{y\thicksim x,y\in\Omega^-}w_{xy}u_n^-(y)u_n^+(x)\right.\\
&   & \left.+\frac{b}{2}\|u_n^\mp\|^2\sum_{x\in\Omega^-}\sum\limits_{y\thicksim x,y\in\Omega^+}w_{xy}u_n^-(x)u_n^+(y)
         +\frac{b}{2}\|u_n^\mp\|^2\sum_{x\in\Omega^+}\sum\limits_{y\thicksim x,y\in\Omega^-}w_{xy}u_n^-(y)u_n^+(x)\right]\\
& = &\int_{\Omega}Q(x)|\tilde{u}^\pm|^{p}\ln |\tilde{u}^\pm|^rd\mu
        +\frac{a}{2}\sum_{x\in\Omega^-}\sum\limits_{y\thicksim x,y\in\Omega^+}w_{xy}\tilde{u}^-(x)\tilde{u}^+(y)
        +\frac{a}{2}\sum_{x\in\Omega^+}\sum\limits_{y\thicksim x,y\in\Omega^-}w_{xy}\tilde{u}^-(y)\tilde{u}^+(x)\\
&    & -\frac{b}{2}\left(\sum_{x\in\Omega^-}
             \sum\limits_{y\thicksim x,y\in\Omega^+}w_{xy}\tilde{u}^-(x)\tilde{u}^+(y)\right)^2
          -\frac{b}{2}\left(\sum_{x\in\Omega^+}
             \sum\limits_{y\thicksim x,y\in\Omega^-}w_{xy}\tilde{u}^-(y)\tilde{u}^+(x)\right)^2\\
&    & -b\|\tilde{u}^+\|^2\|\tilde{u}^-\|^2
         -b\sum_{x\in\Omega^+}\sum\limits_{y\thicksim x,y\in\Omega^-}w_{xy}\tilde{u}^-(y)\tilde{u}^+(x)
             \sum_{x\in\Omega^-}\sum\limits_{y\thicksim x,y\in\Omega^+}w_{xy}\tilde{u}^-(x)\tilde{u}^+(y)\\
&    & +\frac{3b}{2}\|\tilde{u}^\pm\|^2\sum_{x\in\Omega^-}
             \sum\limits_{y\thicksim x,y\in\Omega^+}w_{xy}\tilde{u}^-(x)\tilde{u}^+(y)
         +\frac{3b}{2}\|\tilde{u}^\pm\|^2\sum_{x\in\Omega^+}
             \sum\limits_{y\thicksim x,y\in\Omega^-}w_{xy}\tilde{u}^-(y)\tilde{u}^+(x)\\
&    & +\frac{b}{2}\|\tilde{u}^\mp\|^2\sum_{x\in\Omega^-}
              \sum\limits_{y\thicksim x,y\in\Omega^+}w_{xy}\tilde{u}^-(x)\tilde{u}^+(y)
         +\frac{b}{2}\|\tilde{u}^\mp\|^2\sum_{x\in\Omega^+}
              \sum\limits_{y\thicksim x,y\in\Omega^-}w_{xy}\tilde{u}^-(y)\tilde{u}^+(x),
\end{eqnarray*}
which implies that
\begin{eqnarray*}
\label{aaa1}
\langle I'(\tilde{u}),\tilde{u}^\pm\rangle= 0.
\end{eqnarray*}
Hence, $\tilde{u}\in\mathcal{M}$. Note that $\tilde{u}^\pm\neq0$. If  let $\tilde{s}=0$ and $\tilde{t}=0$ in (\ref{Eq4}), then we have
\begin{eqnarray*}
     m=I(\tilde{u})
\geq a\left(\frac{1}{2}-\frac{1}{p}\right)\|\tilde{u}^+\|^2+a\left(\frac{1}{2}-\frac{1}{p}\right)\|\tilde{u}^-\|^2
     +b\left(\frac{1}{4}-\frac{1}{p}\right)\|\tilde{u}^+\|^4+b\left(\frac{1}{4}-\frac{1}{p}\right)\|\tilde{u}^-\|^4>0.
\end{eqnarray*}
By a similar argument as above, we can also obtain that $c>0$ is achieved.
\qed

\vskip2mm
\noindent
{\bf Lemma 3.9.} {\it The minimizers of  $\inf_{\mathcal{M}}I$ and $\inf_{\mathcal{N}}I$ are critical points of $I$.}
\vskip2mm
\noindent
{\bf Proof.} We prove this conclusion by contradiction. If we let $\tilde{u}=\tilde{u}^++\tilde{u}^-\in\mathcal{M}$, $I'(\tilde{u})=m$ and $I'(\tilde{u})\neq0$,  then there exist $\delta>0$ and $\zeta>0$ such that
$$
\|I'(\tilde{u})\|\geq\zeta,\;\;\;\; \forall \|u-\tilde{u}\|\leq 3\delta\;\;\text{and}\;\;u\in{\mathcal{H}^{1,2}_0(\Omega)}.
$$
Since $\tilde{u}\in\mathcal{M}$, according to the definition of $\mathcal{M}$ and Lemma 3.1, we can obtain that for all $s,t\geq 0$, there exists
\begin{eqnarray}\label{Eq200}
       I(s\tilde{u}^++t\tilde{u}^-)
&\leq& I(\tilde{u})-a\left(\frac{1-s^2}{2}-\frac{1-s^{p}}{p}\right)\|\tilde{u}^+\|^2
         -a\left(\frac{1-t^2}{2}-\frac{1-t^{p}}{p}\right)\|\tilde{u}^-\|^2\nonumber\\
&    & -b\left(\frac{1-s^4}{4}-\frac{1-s^{p}}{p}\right)\|\tilde{u}^+\|^4
        -b\left(\frac{1-t^4}{4}-\frac{1-t^{p}}{p}\right)\|\tilde{u}^-\|^4
        -\frac{b(s^2-t^2)^2}{4}\|\tilde{u}^+\|^2\|\tilde{u}^-\|^2\nonumber\\
&    & -\frac{b(s^2-t^2)^2}{4}\sum_{x\in\Omega^-}\sum\limits_{y\thicksim x,y\in\Omega^+}w_{xy}\tilde{u}^-(x)\tilde{u}^+(y)
           \sum_{x\in\Omega^+}\sum\limits_{y\thicksim x,y\in\Omega^-}w_{xy}\tilde{u}^-(y)\tilde{u}^+(x)\nonumber\\
&    & +\frac{a(s-t)^2}{4}\sum_{x\in\Omega^-}\sum\limits_{y\thicksim x,y\in\Omega^+}w_{xy}\tilde{u}^-(x)\tilde{u}^+(y)
          +\frac{a(s-t)^2}{4}\sum_{x\in\Omega^+}\sum\limits_{y\thicksim x,y\in\Omega^-}w_{xy}\tilde{u}^-(y)\tilde{u}^+(x)\nonumber\\
&    & -\frac{b(s^2-t^2)^2}{8}\left(\sum_{x\in\Omega^-}\sum\limits_{y\thicksim x,y\in\Omega^+}w_{xy}\tilde{u}^-(x)\tilde{u}^+(y)\right)^2\nonumber\\
&    & -\frac{b(s^2-t^2)^2}{8}\left(\sum_{x\in\Omega^+}\sum\limits_{y\thicksim x,y\in\Omega^-}w_{xy}\tilde{u}^-(y)\tilde{u}^+(x)\right)^2\nonumber\\
&\leq& m-a\left(\frac{1-s^2}{2}-\frac{1-s^{p}}{p}\right)\|\tilde{u}^+\|^2
            -a\left(\frac{1-t^2}{2}-\frac{1-t^{p}}{p}\right)\|\tilde{u}^-\|^2\nonumber\\
&    & -b\left(\frac{1-s^4}{4}-\frac{1-s^{p}}{p}\right)\|\tilde{u}^+\|^4
        -b\left(\frac{1-t^4}{4}-\frac{1-t^{p}}{p}\right)\|\tilde{u}^-\|^4
        -\frac{b(s^2-t^2)^2}{4}\|\tilde{u}^+\|^2\|\tilde{u}^-\|^2.\nonumber\\
\end{eqnarray}
Using a discussion method similar to (\ref{G4}) and (\ref{G5}), we can assume the existence of two constants $C_1\in(0,1)$ and $C_2\in(1,+\infty)$ such that
\begin{eqnarray}
\label{Eq21}
\langle I'(s\tilde{u}^++C_1\tilde{u}^-),C_1\tilde{u}^-\rangle>0,\;\;\langle I'(s\tilde{u}^++C_2\tilde{u}^-),C_2\tilde{u}^-\rangle<0,\;\;\forall\;s\in[C_1,C_2],
\end{eqnarray}
and
\begin{eqnarray}
\label{Eq22}
\langle I'(C_1\tilde{u}^++t\tilde{u}^-),C_1\tilde{u}^+\rangle>0,\;\;\langle I'(C_2\tilde{u}^++t\tilde{u}^-),C_2\tilde{u}^+\rangle<0,\;\;\forall\;t\in[C_1,C_2].
\end{eqnarray}
\par
Let $D=\big(C_1,C_2\big)\times\big(C_1,C_2\big)$. Thus, in virtue of (\ref{Eq200}) that
\begin{eqnarray}
\label{Eq23}
\mathcal{X}:=\max_{(s,t)\in\partial D}I(s\tilde{u}^++t\tilde{u}^-)<m.
\end{eqnarray}
For $\varepsilon:=\min\{\frac{m-\mathcal{X}}{3},\frac{\zeta\delta}{8}\}$ and $S_\delta:=B(\tilde{u},\delta)$, by \cite{Minimax}, we can get a deformation $\eta\in\mathcal{C}([0,1]\times \mathcal{H}^{1,2}_0(\Omega),\mathcal{H}^{1,2}_0(\Omega))$ such that\\
$(\text i) \eta(1,u)=u\;\;\text{if}\;\;|I(u)-m|>2\varepsilon$;\\
$(\text {ii}) \eta(1,I^{m+\varepsilon}\cap S_\delta)\subset I^{m-\varepsilon}$,\;\;\text{where}\;\; $I^{c}:=\{u\in \mathcal{H}^{1,2}_0(\Omega):I(u)\leq c\}$;\\
$(\text {iii}) I(\eta(1,u))\leq I(u),\;\;\forall u\in{\mathcal{H}^{1,2}_0(\Omega)}$.\\
It follows from Lemma 3.1, (iii) and for $s,t\geq0$ that make $|s-1|^2+|t-1|^2\geq\delta^2/{\|\tilde{u}\|^2}$ hold,
\begin{eqnarray}
\label{Eq24}
     I(\eta(1,s\tilde{u}^++t\tilde{u}^-))
\leq I(s\tilde{u}^++t\tilde{u}^-)
  <  I(\tilde{u})=m.
\end{eqnarray}
Thus, Corollary 3.3 implies  that $I(s\tilde{u}^++t\tilde{u}^-)\leq I(\tilde{u})=m,\forall s,t\geq0$. Then by (ii), we have
\begin{eqnarray}
\label{Eq25}
I(\eta(1,s\tilde{u}^++t\tilde{u}^-))\leq m-\varepsilon,\;\;\text{for}\;\;s,t\geq0,|s-1|^2+|t-1|^2<\delta^2/{\|\tilde{u}\|^2}.
\end{eqnarray}
According to (\ref{Eq23}), (\ref{Eq24}) and (\ref{Eq25}), we can get that
\begin{eqnarray}
\label{Eq26}
\max_{(s,t)\in\bar{D}}I(\eta(1,s\tilde{u}^++t\tilde{u}^-))<m.
\end{eqnarray}
Define $k(s,t)=s\tilde{u}^++t\tilde{u}^-$. Next, we prove that $\eta(1,k(D))\cap\mathcal{M}\neq\varnothing$. Set $\gamma(s,t):=\eta(1,k(s,t))$ and
\begin{eqnarray*}
      \phi(s,t)
  &:=&  \Big(\left\langle I'(\gamma(s,t)),\gamma^+(s,t)\right\rangle\Big)\in\mathbb{R}^2,\\
  &:=&  (\phi_1(s,t),\phi_2(s,t))\in\mathbb{R}^2,\;\;\forall\;s,t\geq0.
\end{eqnarray*}
In view of (\ref{Eq24}) and the selection of $\varepsilon$, we have $I(s\bar{u}^++t\bar{u}^-)\leq\mathcal{X}\leq m-3\varepsilon<m-2\varepsilon,\forall\;(s,t)\in\partial D$. As a consequence, in virtue of (i), there holds
\begin{eqnarray}
\label{Eq27}
\gamma(s,t)=\eta(1,s\tilde{u}^++t\tilde{u}^-)=s\tilde{u}^++t\tilde{u}^-,\;\;\forall\;(s,t)\in\partial D.
\end{eqnarray}
It follows from (\ref{Eq21}), (\ref{Eq22}) and (\ref{Eq27}) that
\begin{eqnarray}
\label{Eq28}
\phi_1(s,C_1)>0,\;\;\phi_1(s,C_2)<0,\;\;\forall\;s\in[C_1,C_2],
\end{eqnarray}
and
\begin{eqnarray}
\label{Eq29}
\phi_2(C_1,t)>0,\;\;\phi_2(C_2,t)<0,\;\;\forall\;t\in[C_1,C_2].
\end{eqnarray}
Consequently, combining (\ref{Eq28}), (\ref{Eq29}) and the variant of Miranda Theorem in Lemma 2.4 of \cite{Shuai 2017}, there appears a pair of positive numbers $(\bar{s},\bar{t})\in D$ such that $\phi(\bar{s},\bar{t})=0$, that is,
\begin{eqnarray}
\label{Eq30}
\Big(\left\langle I'(\gamma(\bar{s},\bar{t})),\gamma^+(\bar{s},\bar{t})\right\rangle,\left\langle I'(\gamma(\bar{s},\bar{t})),\gamma^-(\bar{s},\bar{t})\right\rangle\Big)=(0,0).
\end{eqnarray}
\par
Note that $\bar{u}\in\mathcal{M}$ and $(\bar{s},\bar{t})\in D$, according to (\ref{Eq30}), we conclude that $\gamma(s,t)=\eta(1,\bar{s}\tilde{u}^++\bar{t}\tilde{u}^-)\in\mathcal{M}$. Thus, $\eta(1,k(D))\cap\mathcal{M}\neq\varnothing$, which contradicts with (\ref{Eq26}) and furthermore, it shows that $I'(\tilde{u})=0$.
\par
Similar to the argument as above, we can  also obtain that the minimizer of $\inf_{\mathcal{N}}I$ is a critical point of $I$.
\qed

\vskip2mm
\noindent
{\bf Proof of Theorem 1.1.} In view of Lemma 3.8 and Lemma 3.9, there exists $\tilde{u}\in\mathcal{M}$ such that $I(\tilde{u})=m$ and $I'(\tilde{u})=0$, and there exists $\bar{u}\in\mathcal{N}$ such that $I(\bar{u})=c$ and $I'(\bar{u})=0$.
Furthermore, according to Corollary 3.3, Proposition 3.2 and Lemma 3.5, there exist positive constants $s_0$ and $t_0$ such that
$s_0\tilde{u}^+\in \mathcal{N}$ and $t_0\tilde{u}^-\in \mathcal{N}$, and then
\begin{eqnarray*}
        m
&  =   &  I(\tilde{u})=\max_{s,t\geq0} I(s\tilde{u}^++t\tilde{u}^-)\\
& \geq &  \max_{s,t\geq0} \big[I(s\tilde{u}^+)+I(t\tilde{u}^-)\big]\\
&  =   &  \max_{s\geq0} I(s\tilde{u}^+)+\max_{t\geq0} I(t\tilde{u}^-)\\
& \geq &  I(s_0\tilde{u}^+)+I(t_0\tilde{u}^-)\\
& \geq &  I(\bar{u})+I(\bar{u})\\
& \geq &  2c>0.
\end{eqnarray*}
The proof of Theorem 1.1 is completed.
\qed

\vskip2mm
{\section{Appendix}}

\noindent{\bf Appendix 1.} {\it $r(1-\tau^p)+p\tau^p\ln \tau^r>0$ for all $\tau\in(0,1)\cup(1+\infty)$ and $p>2$.}
\vskip0mm
\noindent
{\bf Proof.} Let $f(\tau)=r(1-\tau^p)+p\tau^p\ln \tau^r$. Then
$$
    f'(\tau)
 = -rp\tau^{p-1}+p^2\tau^{p-1}\ln \tau^r+rp\tau^{p-1}
 = p^2\tau^{p-1}\ln \tau^r.$$
Take $f'(\tau)=0$. Then we have $\tau=1$, i.e. the function $f(\tau)$ reaches its minimum value at $\tau=1$ and $f(1)=r(1-1)+p1^p\ln 1^r=0$. So,  $r(1-\tau^p)+l\tau^p\ln \tau^r>0$, for all $\tau\in(0,1)\cup(1+\infty)$ and $p>2$.
\qed
\vskip3mm
\noindent{\bf Appendix 2.} {\it The function $f(x)=\frac{1-a^x}{x}$ is strictly monotonically decreasing on $(0,+\infty)$ for $a>0$ and $a\neq 1$.}
\vskip0mm
\noindent
{\bf Proof.} Note that
$$
f'(x)=\frac{-a^x\cdot\ln a\cdot x-(1-a^x)}{x^2}=\frac{a^x-a^x \cdot x\ln a-1}{x^2}.
$$
Obviously, the symbol of $f'(x)$ is the same as that of $k(x):=(a^x-a^x \cdot x\ln a-1)$. Note that $k'(x)=-a^x\cdot x(\ln a)^2<0$ for all $x\in (0,+\infty)$, and so  $k(x)<k(0)=0$ for all $x\in (0,+\infty)$. Then $f'(x)$ monotonically decreasing on $(0,\infty)$ and $f'(0)=\lim_{x\rightarrow0}\frac{a^x-a^x \cdot x\ln a-1}{x^2}=-\frac{(\ln a)^2}{2}<0$ which is the maximum value of $f'(x)$. Thus, $f(x)=\frac{1-a^x}{x}$ is strictly monotonically decreasing on $(0,+\infty)$ for $a>0$ and $a\neq 1$.
\qed
\vskip3mm
\noindent{\bf Funding information}\\
\noindent
This project is supported by Yunnan Fundamental Research Projects (grant No: 202301AT070465) of China and  Xingdian Talent Support Program for Young Talents of Yunnan Province of China.

\vskip3mm
 \noindent
\noindent{\bf Author contribution}
\par
\noindent
 These authors contributed equally to this work.

\vskip3mm
 \noindent
\noindent{\bf Conflict of interest}

\noindent
The authors state no conflict of interest.

\vskip3mm
 \noindent
\noindent{\bf Data availability statement }

\noindent
Not available.

\vskip2mm
\renewcommand\refname{References}
{}
\end{document}